\documentclass[a4paper, 11pt, reqno]{amsart}
\usepackage{a4wide}
\usepackage{amssymb, amsmath, amsthm}

\usepackage{hyperref}

\newtheorem{de}{Definition}[section]
\newtheorem{lem}[de]{Lemma}
\newtheorem{prop}[de]{Proposition}

\newtheorem{thm}[de]{Theorem}

\numberwithin{equation}{section}

\newcommand{\BlackBox}{\rule{1.5ex}{1.5ex}}
\newenvironment{prf}{\par\noindent{\bf Proof.\ }}{\hfill\rule{1.5ex}{1.5ex}\vspace{0.3cm}}

\newcommand{\bprf}{\begin{prf}}
\newcommand{\eprf}{\end{prf}}

\usepackage{amssymb}
\usepackage{dsfont}
\usepackage[T1]{fontenc}

\newcommand{\R}{\mathbb{R}}

\newcommand{\C}{\mathbb{C}}
\newcommand{\N}{\mathbb{N}}

\newcommand{\E}{\mathbb{E}}
\newcommand{\F}{\mathbb{F}}

\newcommand{\Y}{\mathbb{Y}}

\newcommand{\D}{\,\textrm{d}}
\newcommand{\g}{\textrm{g}}

\newcommand{\bea}{\begin{eqnarray*}}
\newcommand{\eea}{\end{eqnarray*}}

\newcommand{\beq}{\begin{equation}}
\newcommand{\eeq}{\end{equation}}

\newcommand{\tr}{\text{tr}}

\newcommand{\ra}{\rightarrow}

\newcommand{\hra}{\hookrightarrow}

\newcommand{\bc}{\begin{mathcal}}
\newcommand{\ec}{\end{mathcal}}
\newcommand{\calF}{\bc F\ec}
\newcommand{\calG}{\bc G\ec}
\newcommand{\calC}{\bc C\ec}
\newcommand{\calP}{\bc P\ec}
\newcommand{\calM}{\bc M\ec}
\newcommand{\calE}{\bc E\ec}
\newcommand{\calL}{\bc L\ec}
\newcommand{\calH}{\bc H\ec}
\newcommand{\calN}{\bc N\ec}
\newcommand{\calA}{\bc A\ec}
\newcommand{\calS}{\bc S\ec}
\newcommand{\calO}{\bc O\ec}
\newcommand{\calT}{\bc T\ec}
\newcommand{\calB}{\bc B\ec}
\newcommand{\calD}{\bc D\ec}
\newcommand{\calR}{\bc R\ec}

\newcommand{\calJ}{\bc J\ec}
\newcommand{\calV}{\bc V\ec}

\newcommand{\HT}{\calH\calT}

\newcommand{\ol}{\overline}
\newcommand{\vphi}{\varphi}

\newcommand{\wt}{\widetilde}
\newcommand{\bs}{\backslash}
\newcommand{\ii}{\textrm{i}}

\newcommand{\supp}{\textrm{supp}}

\newcommand{\rea}{\textrm{Re}\,}

\newcommand{\eps}{\varepsilon}

\def\typeout#1{\message{^^J}\message{#1}\message{^^J}}
\newif\ifSRCOK \SRCOKtrue
\newcount\PAGETOP
\newcount\LASTLINE
\global\PAGETOP=1
\global\LASTLINE=-1
\def\EJECT{\SRC\eject}
\def\WinEdt#1{\typeout{:#1}}
\gdef\MainFile{\jobname.tex}
\gdef\CurrentInput{\MainFile}
\newcount\INPSP
\global\INPSP=0
\def\SRC{\ifSRCOK%
  \ifnum\inputlineno>\LASTLINE%
    \ifnum\LASTLINE<0%
      \global\PAGETOP=\inputlineno%
    \fi%
    \global\LASTLINE=\inputlineno%
    \ifnum\INPSP=0%
      \ifnum\inputlineno>\PAGETOP%
        
      \fi%
    \else%
      
    \fi%
  \fi%
\fi}
\def\PUSH#1{%
\SRC%
\ifnum\INPSP=0 \global\let\INPSTACKA=\CurrentInput \else%
\ifnum\INPSP=1 \global\let\INPSTACKB=\CurrentInput \else%
\ifnum\INPSP=2 \global\let\INPSTACKC=\CurrentInput \else%
\ifnum\INPSP=3 \global\let\INPSTACKD=\CurrentInput \else%
\ifnum\INPSP=4 \global\let\INPSTACKE=\CurrentInput \else%
\ifnum\INPSP=5 \global\let\INPSTACKF=\CurrentInput \else%
               \global\let\INPSTACKX=\CurrentInput \fi\fi\fi\fi\fi\fi%
\gdef\CurrentInput{#1}%
\WinEdt{<+ \CurrentInput}%
\global\LASTLINE=0%
\ifSRCOK\fi%
\global\advance\INPSP by 1}
\def\POP{%
\ifnum\INPSP>0 \global\advance\INPSP by -1  \fi%
\ifnum\INPSP=0 \global\let\CurrentInput=\INPSTACKA \else%
\ifnum\INPSP=1 \global\let\CurrentInput=\INPSTACKB \else%
\ifnum\INPSP=2 \global\let\CurrentInput=\INPSTACKC \else%
\ifnum\INPSP=3 \global\let\CurrentInput=\INPSTACKD \else%
\ifnum\INPSP=4 \global\let\CurrentInput=\INPSTACKE \else%
\ifnum\INPSP=5 \global\let\CurrentInput=\INPSTACKF \else%
               \global\let\CurrentInput=\INPSTACKX \fi\fi\fi\fi\fi\fi%
\WinEdt{<-}%
\global\LASTLINE=\inputlineno%
\global\advance\LASTLINE by -1%
\SRC}
\def\INPUT#1{\relax}

\let\originalxxxeverypar\everypar
\newtoks\everypar
\originalxxxeverypar{\the\everypar\expandafter\SRC}
\everymath\expandafter{\the\everymath\expandafter\SRC}
\output\expandafter{\expandafter\SRCOKfalse\the\output}

\newif\ifSRCOK \SRCOKtrue
\newcount\PAGETOP
\newcount\LASTLINE
\global\PAGETOP=1
\global\LASTLINE=-1
\gdef\MainFile{\jobname.tex}
\gdef\CurrentInput{\MainFile}
\newcount\INPSP
\global\INPSP=0
\def\EJECT{\SRC\eject}
\def\WinEdt#1{\typeout{:#1}}
\def\SRC{\ifSRCOK%
  \ifnum\inputlineno>\LASTLINE%
    \ifnum\LASTLINE<0%
      \global\PAGETOP=\inputlineno%
    \fi%
    \global\LASTLINE=\inputlineno%
    \ifnum\INPSP=0%
      \ifnum\inputlineno>\PAGETOP%
      \fi%
    \else%
    \fi%
  \fi%
\fi}
\def\PUSH#1{%
\SRC%
\ifnum\INPSP=0 \global\let\INPSTACKA=\CurrentInput \else%
\ifnum\INPSP=1 \global\let\INPSTACKB=\CurrentInput \else%
\ifnum\INPSP=2 \global\let\INPSTACKC=\CurrentInput \else%
\ifnum\INPSP=3 \global\let\INPSTACKD=\CurrentInput \else%
\ifnum\INPSP=4 \global\let\INPSTACKE=\CurrentInput \else%
\ifnum\INPSP=5 \global\let\INPSTACKF=\CurrentInput \else%
               \global\let\INPSTACKX=\CurrentInput \fi\fi\fi\fi\fi\fi%
\gdef\CurrentInput{#1}%
\WinEdt{<+ \CurrentInput}%
\global\LASTLINE=0%
\ifSRCOK\fi%
\global\advance\INPSP by 1}
\def\POP{%
\ifnum\INPSP>0 \global\advance\INPSP by -1  \fi%
\ifnum\INPSP=0 \global\let\CurrentInput=\INPSTACKA \else%
\ifnum\INPSP=1 \global\let\CurrentInput=\INPSTACKB \else%
\ifnum\INPSP=2 \global\let\CurrentInput=\INPSTACKC \else%
\ifnum\INPSP=3 \global\let\CurrentInput=\INPSTACKD \else%
\ifnum\INPSP=4 \global\let\CurrentInput=\INPSTACKE \else%
\ifnum\INPSP=5 \global\let\CurrentInput=\INPSTACKF \else%
               \global\let\CurrentInput=\INPSTACKX \fi\fi\fi\fi\fi\fi%
\WinEdt{<-}%
\global\LASTLINE=\inputlineno%
\global\advance\LASTLINE by -1%
\SRC}
\def\INPUT#1{\relax}
\let\OldINCLUDE=\include
\def\include#1{
\EJECT%
\PUSH{#1.tex}%
\OldINCLUDE{#1}%
\POP}

\let\originalxxxeverypar\everypar
\newtoks\everypar
\originalxxxeverypar{\the\everypar\expandafter\SRC}
\everymath\expandafter{\the\everymath\expandafter\SRC}
\let\zzzxxxbibliography=\bibliography
\def\bibliography#1{\PUSH{\jobname.bbl}\zzzxxxbibliography{#1}\POP}
\output\expandafter{\expandafter\SRCOKfalse\the\output}

\begin{document}

\title{Maximal regularity with temporal weights for parabolic problems with 
inhomogeneous boundary conditions}

\author{Martin Meyries}
\author{Roland Schnaubelt}

\address{Department of Mathematics,
Karlsruhe Institute of Technology, 76128 Karlsruhe, Germany.}
\email{martin.meyries@kit.edu}
\email{schnaubelt@kit.edu}

\thanks{This paper is part of a research project supported by  the Deutsche 
Forschungsgemeinschaft (DFG)}




\begin{abstract}
We develop a maximal regularity approach in temporally weighted
 $L_p$-spaces for vector-valued parabolic initial-boundary value problems with inhomogeneous 
 boundary conditions, both of static and of relaxation type. Normal ellipticity and conditions of 
 Lopatinskii-Shapiro type are the basic structural assumptions. The weighted framework  allows
 to reduce the initial regularity and to avoid compatibility conditions at the boundary, and 
 it provides an inherent smoothing effect of the solutions.  Our main tools are 
 interpolation and trace theory for anisotropic Slobodetskii spaces with temporal weights,
 operator-valued functional calculus, as well as localization and perturbation arguments.
 \end{abstract}

\maketitle

\section{Introduction}
In recent years parabolic equations with fully nonlinear boundary conditions
have attracted a lot of interest since they arise in the analysis of free
boundary value problems such as the Stefan problem with surface tension, see e.g.\
\cite{EPS03}, \cite{KPW11} and  \cite{PSZ11}.  These papers use an $L_p$-approach to such problems
which yields strong solutions in maximal regularity classes. In this framework the boundary 
conditions are attained in a classical sense up to initial time, and not just weakly.
This approach is based on linearization and on a sharp $L_p$-regularity theory for linear
inhomogeneous initial-boundary value problems, as  established in \cite{DHP03}, \cite {DHP07}
and \cite{DPZ08} by Denk, Hieber, Pr\"uss and Zacher. Besides the usual static boundary 
conditions, one also has to treat dynamical boundary conditions of relaxation type
which arise in the context of the Stefan problem with surface tension and in related problems.

However, this approach requires  regularities  of the initial values (and hence of the 
nonlinear phase spaces) which are stronger than the norms one can control by standard 
a priori estimates for the nonlinear problems. In related situations
it is known that one can reduce the required initial regularity by means of temporal weights.
In the $L_p$-setting, it is natural to work  in 
\beq\label{p1}
L_{p,\mu}(J;X) := \big\{u:J\ra X \; :\; t^{1-\mu}u\in L_p(J;X)\big\}
\eeq
endowed with its natural norm, where $p\in(1,\infty)$, $\mu\in(1/p,1]$, $T\in(0,\infty]$, 
$J:=(0,T)$, and $t^{1-\mu}u$ denotes the function $t\mapsto t^{1-\mu} u(t)$ on $J$.
The corresponding weighted Sobolev spaces are defined by 
\beq\label{p2}
W_{p,\mu}^k(J;X) := \big \{u:J\ra X \; : \; u,u',..., u^{(k)} \in L_{p,\mu}(J;X)\big\}
\eeq
for $k\in \N$. These spaces and the corresponding anisotropic spaces like $\E_{u,\mu}(J)$ 
defined below are studied by the authors in detail in \cite{MS11}. To a large extent they enjoy 
analogous properties as the corresponding unweighted spaces. 

To see the effect of the weight, we consider a generator $-A$ of an analytic semigroup 
on a Banach space $X$. Then the orbit $u(t)=e^{-tA}u_0$ belongs to the `maximal regularity space'
\[W_{p,\mu}^1(\R_+; X)\cap L_{p,\mu}(\R_+; D(A))\]
if and only if the initial value $u_0$ belongs to the real interpolation space
\[(X,D(A))_{\mu-1/p,p}\,,\]
see e.g.\ Theorem~1.14.5 in \cite{Tri94}. Recall that  one often fixes a large $p\in (1,\infty)$ 
to treat nonlinearities. Hence, in the unweighted case $\mu=1$ the resulting 
initial regularity  is close to $D(A)$. On the other hand,  taking $\mu$ near $1/p$
one almost reaches the base space $X$. Further, for Banach spaces $X$ of class $\calH\calT$ 
(see Section~2), Pr\"uss and Simonett have proved in \cite{PS04} that the inhomogeneous 
evolution equation
\beq\label{I1}
u'(t) + Au(t) = f(t), \quad t>0, \qquad u(0) = 0,
\eeq
has a unique solution in $W_{p,\mu}^1(\R_+; X)\cap L_{p,\mu}(\R_+; D(A))$ for each 
$f\in L_{p,\mu}(J; X)$ if and only if this fact holds for the unweighted case $\mu=1$.
Since the unweighted case is well understood, see e.g.\ \cite{KW04}, 
the $L_{p,\mu}$-approach is quite convenient for parabolic problems covered by \eqref{I1}.

Unfortunately, it seems that a sharp regularity theory for inhomogenous boundary value problems 
is not possible within the abstract framework of an evolution equation like \eqref{I1}. Instead one
has to restrict to a PDE setting. So  we investigate vector-valued linear parabolic 
systems with inhomogeneous boundary conditions, such of static type, i.e.,
\begin{alignat}{3}
\partial_t u + \calA(t,x,D) u & =    f(t,x), & \qquad & x\in  \Omega, & \qquad t\in J,&    
\nonumber\\
\calB_j(t,x,D) u & =  g_j(t,x), && x\in \Gamma,  & t\in J,&  \qquad j=1,...,m, \label{statmrI}\\
u(0,x) & =  u_0(x), &&  x\in \Omega, &&  \nonumber
\end{alignat} 
as well as such of relaxation (or dynamic) type, i.e.,
\begin{alignat}{3}
\partial_t u + \calA(t,x,D) u & =    f(t,x), & \qquad & x\in \Omega, & \qquad &  t\in J,  \nonumber\\
\partial_t \rho + \calB_0(t,x,D) u + \calC_0(t,x,D_\Gamma) \rho & =  g_0(t,x), &&x\in  \Gamma,  &&  t\in J,\quad   \nonumber\\
\calB_j(t,x,D) u + \calC_j (t,x,D_\Gamma)\rho & =  g_j(t,x), &&x\in \Gamma,  &&  t\in J,\qquad j=1,...,m,   \label{dynmaxregI}\\
u(0,x) & =  u_0(x), &&  x\in \Omega, &&  \nonumber \\
\rho(0,x) & =  \rho_0(x), &&  x\in \Gamma. &&  \nonumber
\end{alignat}
It is assumed that $\Omega\subset \R^n$ is a (inner or outer) domain with compact smooth boundary 
$\Gamma = \partial\Omega$. In (\ref{statmrI}) and (\ref{dynmaxregI}) the unknown $u=u(t,x)$ 
takes values in a Banach space $E$, and in (\ref{dynmaxregI}) the additional unknown 
$\rho = \rho(t,x)$, which only lives on the boundary $\Gamma$, takes values in another Banach 
space $F$. Throughout we assume that $E$ and $F$ are of class $\calH\calT$; for instance
$E$ and $F$ can be finite dimensional leading to usual parabolic systems. The differential
 operator $\calA$ is of order $2m$, where $m\in \N$, and $\calB_j$ are corresponding boundary 
 operators of order $m_j$ not larger than $2m-1$. In (\ref{dynmaxregI}) the differential 
 operators $\calC_j$ contain tangential derivatives of any order up to $k_j\in \N_0$. We assume 
 certain ellipticity and Lopatinskii-Shapiro type conditions and  impose regularity 
 conditions on the  coefficients that are appropriate for the applications to quasilinear problems, 
 see e.g.\ \cite{Mey10}, \cite{Mey11}. The details are described  in Section 2. 
 
 We look for strong solutions $u$ of (\ref{statmrI}), resp.\ $(u,\rho)$ of
 (\ref{dynmaxregI}), which satisfy the respective equations pointwise almost everywhere. 
 In particular, $u$ shall belong to 
$$\E_{u,\mu}(J) := W_{p,\mu}^1\big (J; L_p(\Omega;E)\big)\cap L_{p,\mu}\big(J; 
   W_p^{2m}(\Omega;E)\big).$$ 
The space $\E_{\rho,\mu}(J)$ for $\rho$ is chosen in accordance to the structure of 
(\ref{dynmaxregI}) and to the trace theorems established in our paper \cite{MS11}, see Section~2.
In our main results Theorems~\ref{statmrthm} and \ref{dynmrthm} we show the maximal 
$L_{p,\mu}$-regularity for (\ref{statmrI}) and (\ref{dynmaxregI}) on $J$. This means that there 
are data spaces $\calD_{\text{stat}}(J)$ and $\calD_{\text{rel}}(J)$ such that  
(\ref{statmrI})  and (\ref{dynmaxregI}) have  unique solutions $u\in \E_{u,\mu}(J)$ and
$(u,\rho)\in \E_{u,\mu}(J)\times \E_{\rho,\mu}(J)$, respectively,
if and only if the data satisfies
\[(f,g_1, ..., g_m, u_0)\in \calD_{\text{stat}}(J) \quad\text{and} \quad
(f,g_0,g_1,..., g_m, u_0,\rho_0)\in \calD_{\text{rel}}(J),\]
respectively. The data spaces contain the necessary regularities for the data
and their compatibility conditions at $t=0$ enforced by the static and dynamical
boundary equations in  (\ref{statmrI})  and (\ref{dynmaxregI}). The precise formulations
of these spaces is suggested by the space-time trace theorems from \cite{MS11}.
In the unweighted case $\mu=1$ and with essentially the same assumptions, the maximal 
$L_p$-regularity for (\ref{statmrI}) and (\ref{dynmaxregI}) has been shown by Denk, Hieber $\&$ 
Pr\"uss \cite{DHP07} and Denk, Pr\"uss $\&$  Zacher \cite{DPZ08}, respectively.  
 
We note that the  problem (\ref{dynmaxregI}) is more involved in the several respects.
Clearly, it contains a second variable and a second evolutionary equation. Moreover, 
the operators $\calC_j$ can make the main parts of the equations at the boundary highly 
non homogeneous which then leads to a rather sophisticated solution space $\E_{\rho,\mu}(J)$ 
and to a complicated analysis. It further can happen that $\partial_t\rho$ is continuous
in $t$ up to $t=0$ so that the dynamical equation for $\rho$ leads to an additional compatibility
condition on the regularity of
$\calB_0(0,\cdot,D) u_0 + \calC_0(0,\cdot,D_\Gamma) \rho_0 - g_0(0,\cdot)$.

The main feature of the weighted approach is the flexibility for the regularity of the initial 
values as $\mu$ varies in $(1/p,1]$. We describe these properties in more detail at the end 
of Section~2, indicating here the basic points. We can solve (\ref{statmrI}) and 
(\ref{dynmaxregI}) with the Besov regularity $u_0\in B_{p,p}^{2(\mu-1/p)}(\Omega;E)$ which  
approaches $L_p(\Omega;E)$ as $\mu$ tends to $1/p$. Moreover, if the initial regularity
is sufficently low we loose the compatibility conditions such as
$\calB_j(0,\cdot,D) u_0 + \calC_j(0,\cdot,D_\Gamma) \rho_0 = g_j(0,\cdot)$.
Since the weight has an influence only at $t=0$, our approach yields an inherent smoothing effect 
for the solutions. In particular, for (\ref{statmrI}) one can control the norm of $u(t)$ in 
$B_{p,p}^{2(1-1/p)}(\Omega;E)$ by the much lower norm of $u_0$ in $B_{p,p}^{2(\mu-1/p)}(\Omega;E)$.
For bounded $\Omega$, this fact gives the important compactness of the semiflow solving the 
related nonlinear problems. Also for unbounded $\Omega$ one can thus `upgrade' the usual a priori 
estimates in low norms up to $B_{p,p}^{2(1-1/p)}(\Omega;E)$ if one is able to handle the 
involved nonlinearities. See \cite{Mey10}, \cite{Mey11} and also 
\cite{KPW10} in the framework of \cite{PS04}, as well as \cite{KPW11}. In these papers the weighted
approach was used to establish convergence to equilibria and the existence of global attractors
in high norms.

In Sections 3 and 4 we prove Theorems \ref{statmrthm} and \ref{dynmrthm}.
We first consider model problems with homogeneous constant coefficient operators on the 
full-space $\R^n$  and on the half-space $\R_+^n$ in Section~3. The full-space case, where 
boundary conditions are not involved, is treated by means of \cite{PS04}. For the half-space case 
with boundary conditions we apply the Fourier transform with respect to time and space. The 
solution operators for the resulting ordinary initial value problems have been analysed in
\cite{DHP07} and \cite{DPZ08} for the unweighted case. We now use a recent operator-valued Fourier 
multiplier theorem in the  $L_{p,\mu}$-spaces due to Girardi $\&$ Weis \cite{GW11} and several
  isomorphisms acting on a  scale of weighted anisotropic fractional order spaces which are 
investigated in  \cite{MS11}. It is then possible  to invert the Fourier transforms and to solve 
the half-space 
 problem in the required norms. The case of a general domain is then a consequence of perturbation 
 and localization arguments, and it is considered in Section 4.
 
 Finally we discuss several important special cases of (\ref{statmrI}) and (\ref{dynmaxregI})
 arising as
linearizations of various quasilinear parabolic problems with nonlinear static or 
dynamic boundary conditions.  For instance, the linearization of a reaction-diffusion system 
with nonlinear Robin boundary conditions is of the form (\ref{statmrI}) for
$$E= \R^N, \qquad m=1, \qquad \calA(D) = -\Delta, \qquad \calB_1(x,D) = \partial_\nu := \nu \cdot \tr_\Omega \nabla,$$
where $\Delta$ is the Laplacian and $\nu$ denotes the outer unit normal field on $\Gamma$. 
The linearization of Cahn-Hilliard phase-field models  leads to similar problems of order $4$ 
(i.e., $m=2$). If we take
$$\calB_1 = \tr_\Omega, \qquad \calC_1 = -1,$$
the static boundary condition for $j=1$ in (\ref{dynmaxregI}) reads $u|_\Gamma=\rho$, which leads 
to inhomogeneous dynamic boundary conditions. Hence $\rho$ is simply the spatial trace of $u$ in 
this case.  Now one can take $\calC_0(x,D_\Gamma) = -\Delta_\Gamma$, the Laplace-Beltrami 
operator on $\Gamma$, to obtain boundary conditions describing surface diffusion, i.e.,
$$\partial_t u|_\Gamma + \partial_\nu u -\Delta_\Gamma  u|_\Gamma  =  g_0 \qquad \text{on } \Gamma, \qquad t\in J.$$
If we choose
$$u|_\Gamma + \Delta_\Gamma \rho = g_1 \qquad \text{on } \Gamma, \qquad t\in J,$$ 
as the first static boundary condition in (\ref{dynmaxregI}), we arrive at the linearization of 
the Stefan problem with surface tension as studied in \cite{EPS03}. Here the graph of 
$\rho(t,\cdot)$ is related to the unknown boundary at time $t$. We refer to Section~3 of 
\cite{DPZ08} for further interesting problems that may be written in the form 
(\ref{dynmaxregI}).

\medskip

\textbf{Notations.} We write $a\lesssim b$ for some quantities $a,b$ if there is a generic positive constant $C$ with $a\leq Cb$. If $A$ is a sectorial operator on a Banach space $E$, $\theta\in (0,1)$ and $q\in [1,\infty]$, then we set $D_L(\theta,p) := \big( E,D(L)\big)_{\theta,q}$ for the real interpolation scale between $E$ and $D(L)$. If $X,Y$ are Banach spaces we denote by $\calB(X,Y)$ the space of bounded linear operators between them, with $\calB(X) := \calB(X,X)$.

\section{The assumptions and the results}
Throughout we assume that the Banach spaces $E,F$ are of class $\calH\calT$ (or, equivalently, 
are UMD spaces). This means that the Hilbert transform is bounded on $L_2(\R;X)$ which holds, 
e.g., in Hilbert spaces $X$ or if $X$
is a reflexive Lebesgue or (fractional) Sobolev space; see Sections~III.4.3-4.5 of \cite{Ama95}.
We first describe the differential 
operators in (\ref{statmrI}) and (\ref{dynmaxregI}) in detail. For both problems the operator 
$\calA$ is given by 
$$\calA(t,x,D) = \sum_{|\alpha|\leq 2m} a_\alpha(t,x) D^\alpha, 
           \qquad x\in \Omega, \qquad t\in J,$$ 
where $m\in \N$ and $D = -\ii\nabla$ with the euclidian gradient 
$\nabla= (\partial_{x_1}, ..., \partial_{x_n})$ on $\R^n$. Hence the order of $\calA$ is $2m$. 
The coefficients take values in the bounded linear operators on $E$, i.e., 
$a_\alpha(t,x)\in \calB(E)$. Also for both problems the boundary operators $\calB_j$ are 
of the form
$$\calB_j(t,x,D) = \sum_{|\beta|\leq m_j} b_{j\beta}(t,x) \tr_\Omega D^\beta, 
   \qquad x\in \Gamma, \qquad t\in J, \qquad j=0,...,m,$$ 
where $m_j\in \{0,..., 2m-1\}$ is the order of $\calB_j$ and the coefficients satisfy
$b_{0\beta}(t,x) \in \calB(E,F)$ and $b_{j\beta}(t,x)\in \calB(E)$ for $j=1,...,m.$ 
We note that  $\calB_j$ acts on $u$ by applying first the euclidian derivatives and then the 
spatial trace $\tr_\Omega$. We assume that each of these operators is nontrivial, i.e., 
$\calB_j\neq 0$ for all $j$.

In problem (\ref{dynmaxregI}), the boundary conditions of relaxation type involve another 
set of operators $\calC_0,..., \calC_m$, which act only on $\rho$ in the following way. 
For  $t\in J$ it is assumed that $\calC_j(t,\cdot, D_\Gamma)$ is a linear map 
$$C^\infty(\Gamma;F)\ra L_1(\Gamma;F)$$
such that for all $j=0,...,m$, all local coordinates $\g$ for $\Gamma$  and all 
$\rho\in C^\infty(\Gamma;F)$ it holds
$$\big (\calC_j(t,\cdot,D_\Gamma)\rho\big ) \circ \g(x) 
   = \sum_{|\gamma|\leq k_j} c_{j\gamma}^{\g}(t,x) D_{n-1}^\gamma (\rho\circ \g)(x), 
     \qquad x\in \g^{-1}(\Gamma\cap U), \qquad t\in J,$$
where $U\subset \R^n$ is the domain of the chart corresponding to $\g$. Here we have 
$D_{n-1} = -\ii\nabla_{n-1}$ with the euclidian gradient $\nabla_{n-1}$ on $\R^{n-1}$.
The order $k_j\in \N_0$ of $\calC_j$ is given independently of the orders of $\calA$ and the $\calB_j$. 
The local coefficients $c_{j\gamma}^{\g}$, that may depend on the  coordinates $\g$, 
are assumed to satisfy
$c_{0\gamma}^{\g}(t,x) \in \calB(F)$ and $c_{j\gamma}^{\g}(t,x)\in \calB(F,E)$ for $j=1,...,m.$
It is assumed that at least one operator $\calC_j$ is nontrivial. If an operator $\calC_j$ 
is trivial, i.e., $\calC_j \equiv 0$, then we set 
$k_j := -\infty$ as its order. Note that we do not assume that an operator $\calC_j$ has global 
coefficients, in the sense that there are functions $c_{j\gamma}$ on $\Gamma$ satisfying 
$c_{j\gamma}^{\g} = c_{j\gamma}\circ \g$ in all coordinates $\g$. In contrast to that, the 
coefficients of the $\calB_j$ are globally defined on $\Gamma$. The standard examples for such 
an operator $\calC_j$ are the Laplace-Beltrami operator $\Delta_\Gamma$ and a convection 
term $\calV \cdot \nabla_\Gamma$, where $\calV$ is a tangential vector field and 
$\nabla_\Gamma$ is the surface gradient on $\Gamma$. Throughout we let
$$p\in (1,\infty) \quad \text{ and }\quad \mu\in (1/p,1].$$ 
We look for solutions $u$ of (\ref{statmrI}) and $(u,\rho)$ of (\ref{dynmaxregI}) 
such that 
$$u\in \E_{u,\mu}= W_{p,\mu}^1\big (J;L_p(\Omega;E)\big) 
                     \cap L_{p,\mu}\big(J; W_p^{2m}(\Omega;E)\big).$$ 
The weighted vector-valued $L_{p,\mu}$-spaces and the corresponding Sobolev spaces spaces 
$W_{p,\mu}^1$ are defined in (\ref{p1}) and (\ref{p2}), respectively, and 
$W_p^{2m}(\Omega;E)$ is the $E$-valued Sobolev space of order $2m$ over $\Omega$. 
For such solutions $u$ the differential equation on the domain $\Omega$ holds for a.e.\ $(t,x)$.
The regularity of $u$  yields
$$f \in \E_{0,\mu} :=L_{p,\mu}(J; L_p(\Omega;E)).$$ 
{}From the mapping properties of the temporal trace described in
Theorem~4.2 of \cite{MS11}, we deduce 
$$u|_{t=0} =u_0 \in X_{u,\mu}:= B_{p,p}^{2m(\mu-1/p)}(\Omega;E).$$
Here $B_{p,p}^{s}(\Omega;E)$ denotes the $E$-valued Besov spaces over $\Omega$. We refer to
 \cite{Ama97}, \cite{Sch} or \cite{Zi89} for a definition and the properties of these spaces. 
Further, Lemma~3.4 of \cite{MS11} shows  that  the operator $D^\beta$ maps $\E_{u,\mu}$ 
continuously into
\begin{equation}\label{diff-map}
 H_{p,\mu}^{1-m_j/2m}\big( J; L_p(\Omega;E)\big)\cap L_{p,\mu}\big(J; W_p^{2m-m_j}(\Omega;E)\big)
 \end{equation}
 for $|\beta|\leq m_j\leq 2m-1$.
 Due to Theorem~4.6 of \cite{MS11}, the spatial trace 
$\tr_\Omega$ maps the space in \eqref{diff-map} continuously into
$$W_{p,\mu}^{\kappa_j}\big (J; L_p(\Gamma;E)\big)\cap 
             L_{p,\mu}\big(J; W_p^{2m\kappa_j}(\Gamma;E)\big), \qquad j=0,...,m,$$
where the number $\kappa_j$ is defined by
$$\kappa_j := 1- \frac{m_j}{2m} -\frac{1}{2mp}, \qquad j=0,...,m.$$
Below we  assume that $\kappa_j \neq 1-\mu + 1/p$ for all $j=0,...,m.$
The  weighted Sobolev spaces $H_{p,\mu}^s$ and the Slobodetskii spaces 
$W_{p,\mu}^s$ of fractional order $s\geq 0$ are  defined by complex and real
interpolation, respectively. The properties of $W_{p,\mu}^s$ and $H_{p,\mu}^s$ are studied in
\cite{MS11}. Moreover, $W_p^s(\Gamma;E)$ is the $E$-valued Sobolev-Slobodetskii space of order $s$, 
where $W_p^s = B_{p,p}^s$ for $s\notin \N_0$. Since the dynamic boundary equation in 
(\ref{dynmaxregI}) takes place in $F$ and the static boundary equations in (\ref{statmrI}) 
and (\ref{dynmaxregI}) take place in $E$, these considerations suggest that we should 
look at boundary data
\begin{align*}
g_0 &\in \F_{0,\mu}:= W_{p,\mu}^{\kappa_0}\big (J; L_p(\Gamma;F)\big)\cap 
          L_{p,\mu}\big(J; W_p^{2m\kappa_0}(\Gamma;F)\big),\\
g_j &\in \F_{j,\mu}:= W_{p,\mu}^{\kappa_j}\big(J; L_p(\Gamma;E)\big)\cap 
        L_{p,\mu}\big(J; W_p^{2m\kappa_j}(\Gamma;E)\big),\qquad j=1,...,m.
\end{align*}        
For convenience we write 
$$\F_\mu := \F_{0,\mu} \times ...\times \F_{m,\mu}, \qquad g= (g_0,..., g_m)\in \F_\mu,$$
and similiarly $\wt{\F}_\mu := \F_{1,\mu} \times ...\times \F_{m,\mu}$ and
$\wt{g}= (g_1,..., g_m)\in \wt{\F}_\mu$.

We now determine the regularity of $\rho$ and $\rho_0$ in (\ref{dynmaxregI}). 
Assuming sufficient smoothness of the coefficients of the operators, we look for
a space $\E_{\rho,\mu}$ for $\rho$ such that each term in (\ref{dynmaxregI}) 
involving $\rho$  acts continuously from $\E_{\rho,\mu}$ to the space $\F_{j,\mu}$ 
where the respective equation takes place. 
It can be seen as in Section~2 of \cite{DPZ08} that
\begin{align}\label{E-rho}
\E_{\rho,\mu}= &\, W_{p,\mu}^{1+\kappa_0}\big (J; L_p(\Gamma;F)\big)\cap 
         L_{p,\mu}\big (J; W_p^{l + 2m\kappa_0}(\Gamma;F)\big ) \vspace{0.2cm}\\
& \quad\cap W_{p,\mu}^1\big (J; W_p^{2m\kappa_0}(\Gamma;F)\big ) \cap 
\sideset{}{_{j\in \wt{\calJ}}}\bigcap W_{p,\mu}^{\kappa_j}\big (J; W_p^{k_j}(\Gamma;F)\big ) \notag
\end{align} 
satisfies these requirements, where we put 
$$\wt{\calJ} := \big \{ j\in \{0,...,m\}\;:\; k_j\neq -\infty\big \}, \qquad l_j := k_j-m_j+m_0, \qquad l:= \max_{j=0,...,m} l_j.$$ 
It is important to note that 
\begin{equation}\label{kappa-l} 
2m\kappa_j+k_j = 2m\kappa_0 +l_j, \qquad j=0,\dots,m.
\end{equation}
We represent
$\E_{\rho,\mu}$ by the points $(0,1+\kappa_0)$, $(l+2m\kappa_0,0)$, $(2m\kappa_0,1)$ and
 $(k_j,\kappa_j)$, $j\in \wt{\calJ}$, corresponding to the space-time differentiability of the 
 spaces $Z_i$ on the right-hand side of \eqref{E-rho}. The \textsl{Newton polygon} $\calN\calP$ for 
 $\E_{\rho,\mu}$ is then defined as the convex hull of these points together with $(0,0)$. 
The \textsl{leading part} of $\calN\calP$  is the polygonal part of its boundary
 connecting $(0,1+\kappa_0)$ to $(l+2m\kappa_0,0)$ anti-clockwise.
 
 Let $Z_i$ and $Z_j$ be two different spaces on  the right-hand side of \eqref{E-rho}. It is 
 shown in Proposition~3.2 of \cite{MS11} that $Z_i\cap Z_j$ embeds into all spaces 
 whose space-time regularity corresponds to the line segment connecting the two points
 that represent $Z_i$ and $Z_j$ in $\calN\calP$. Consequently, the description of  
$\E_{\rho,\mu}$ given in \eqref{E-rho} contains redundant spaces, in general. We derive
a  nonredundant description of $\E_{\rho,\mu}$  as in the case $\mu=1$ presented 
in \cite{DPZ08}. Here one has to distinguish three cases. In each case, a direct application
Theorem~4.2 of \cite{MS11} further yields the temporal trace space of $\rho$ at $t=0$,
 denoted by 
$$X_{\rho,\mu}:= \tr_{t=0} \E_{\rho,\mu}.$$ 
In the same way we also obtain that  the temporal derivative $\partial_t \rho$ has a trace 
at $t=0$ if $\kappa_0> 1-\mu+1/p$. We denote the resulting trace space by
$$X_{\partial_t\rho,\mu}:= \tr_{t=0} \partial_t \E_{\rho,\mu} \qquad 
   \text{if }\; \kappa_0 >1-\mu+1/p.$$ 
We remark that Theorem~4.2 of \cite{MS11} means that these trace spaces are given by 
$B_{p,p}^\sigma(\Gamma;F)$
for the numbers $\sigma>0$ such that $(\sigma,k+1-\mu+1/p)$ belongs to leading part
of $\calN\calP$ for $k=0$ and $k=1$, respectively.
We can now give the nonredundant description of the spaces $\E_{\rho,\mu}$, $X_{\rho,\mu}$ and $X_{\partial_t \rho,\mu}$.\vspace{0.1cm}\\
\noindent \textbf{Case 1: $\boldsymbol{l=2m}$.} One has 
$$\E_{\rho,\mu} = W_{p,\mu}^{1+\kappa_0}\big (J; L_p(\Gamma;F)\big) 
      \cap L_{p,\mu}\big(J; W_p^{2m(1+\kappa_0)}(\Gamma;F)\big)$$ 
since all other spaces in \eqref{E-rho} correspond to points on or below the straight line
$s=1+\kappa_0-r/2m$ from $(0,1+\kappa_0)$ to $(2m+2m\kappa_0,0)$ due to \eqref{kappa-l}.
It follows that
$$X_{\rho,\mu} = B_{p,p}^{2m(\kappa_0 + \mu-1/p)}(\Gamma;F), \qquad X_{\partial_t \rho,\mu} = B_{p,p}^{2m(\kappa_0 -(1- \mu+1/p))}(\Gamma;F)\quad \text{if }\; \kappa_0 >1-\mu+1/p.$$
\textbf{Case 2: $\boldsymbol{l<2m}$.} One has
\begin{align*}
 \E_{\rho,\mu}  = \quad & W_{p,\mu}^{1+\kappa_0}\big(J; L_p(\Gamma;F)\big) \cap  L_{p,\mu}\big(J; W_p^{l+2m\kappa_0}(\Gamma;F)\big) \cap W_{p,\mu}^1\big(J; W_p^{2m\kappa_0}(\Gamma;F)\big)
\end{align*}
since $(1,2m\kappa_0)$ lies above the line segment $s=1+\kappa_0-r(1+\kappa_0)/(l+2m\kappa_0)$
from $(0,1+\kappa_0)$ to $(2m+2m\kappa_0,0)$
and all points $(\kappa_j, k_j)$ are below the line $s=1+(2m\kappa_0-r)/l$ 
connecting $(1,2m\kappa_0)$ and $(0, l+2m\kappa_0)$. We then obtain the trace spaces
$$X_{\rho,\mu} = B_{p,p}^{2m\kappa_0 + l(\mu-1/p)}(\Gamma;F), \qquad 
X_{\partial_t \rho,\mu} = B_{p,p}^{2m(\kappa_0 -(1- \mu+1/p))}(\Gamma;F)
 \quad \text{if }\; \kappa_0 >1-\mu+1/p.$$
\textbf{Case 3: $\boldsymbol{l>2m}$.}  Now $(1,2m\kappa_0)$ belongs to the interior
of $\calN\calP$ and it holds
\begin{align*}
 \E_{\rho,\mu} = \quad & W_{p,\mu}^{1+\kappa_0}\big(J; L_p(\Gamma;F)\big)\cap  
  L_{p,\mu}\big(J; W_p^{l+2m\kappa_0}(\Gamma;F)\big) \cap \sideset{}{_{j\in \calJ}}\bigcap  
 W_{p,\mu}^{\kappa_j}\big(J; W_p^{k_j}(\Gamma;F)\big),
\end{align*}
where $\calJ = \{j_1,..., j_{q_{\text{max}}}\}\subset \wt{\calJ}$, $q_{\text{max}}\in \N$, 
contains those indices $j\in \wt{\calJ}$ so that $(k_j, \kappa_j)$ belongs to the leading part 
of the Newton polygon, i.e., the points 
$$P_0 = (0, 1+\kappa_0), \ \ P_1 = (k_{j_1}, \kappa_{j_1}), \ \ \dots, \ \ 
P_{q_{\text{max}}} = (k_{j_{q_{\text{max}}}}, \kappa_{j_{q_{\text{max}}}}), \ \
P_{q_{\text{max}+1}} = (l+2m\kappa_0, 0)$$
are the vertices of the leading part. It is assumed that $\calJ$  is ordered so that 
$k_{j_{q_1}} < k_{j_{q_2}}$ and $\kappa_{j_{q_1}} > \kappa_{j_{q_2}}$ for $q_1 < q_2.$ 
We also define
$k_{-1}:= 0$, $\kappa_{-1}:=1+\kappa_0$, $m_{-1} := m_0-2m$ and $l_{-1} := 2m.$ 
We further denote the edge in the Newton polygon connecting the points $P_q$ and $P_{q+1}$ by 
$\calN\calP_q$, $q=0,..., q_{\text{max}}$, and set
\begin{align*}
\calJ_{2q}&\,:= \big\{ j\in \calJ \cup \{-1\} \;:\; (k_j,\kappa_j) = P_q\big \}, \qquad q = 0,..., q_{\text{max}},\\
\calJ_{2q+1}&\,:= \big\{ j\in \calJ \cup \{-1\} \;:\; (k_j,\kappa_j) \in  \calN\calP_q\big \}, \qquad q = 0,..., q_{\text{max}}.
\end{align*}
The temporal trace space of $\partial_t\rho$ is obtained by Theorem~4.2 of \cite{MS11} from 
the spaces corresponding to $P_0$ and $P_1$. 
We thus deduce
$$X_{\partial_t \rho,\mu} = B_{p,p}^{k_{j_1}(\kappa_0-(1-\mu+1/p))/(1+\kappa_0 - \kappa_{j_1})}(\Gamma;F)\quad \text{if }\; \kappa_0 >1-\mu+1/p.$$ 
For $X_{\rho,\mu}$ one has to distinguish three more cases.

\noindent \textbf{Case 3(i):} If $\kappa_j>1-\mu+1/p$ for all $j\in \calJ$ then 
$$X_{\rho,\mu} = B_{p,p}^{l + 2m(\kappa_0 - (1-\mu+1/p))}(\Gamma;F).$$
Here we apply Theorem~4.2 of \cite{MS11} to the spaces corresponding to 
$P_{q_{\text{max}}}$ and $P_{q_{\text{max}+1}}.$

\noindent  \textbf{Case 3(ii):} Denote by $j_{q_1}\in \calJ$ be the smallest index with 
$\kappa_{j_{q_1}}> 1-\mu+ 1/p$, and by $j_{q_2}\in \calJ$ the largest index with 
$\kappa_{j_{q_2}}<1-\mu+ 1/p$. Using the spaces corresponding to these indices, 
we conclude that
$$X_{\rho,\mu} = B_{p,p}^{k_{j_{q_1}} + (\kappa_{j_{q_1}} - (1-\mu+1/p)) 
     \frac{k_{j_{q_2}} - k_{j_{q_1}}}{\kappa_{j_{q_2}} - \kappa_{j_{q_1}}}}(\Gamma;F).$$
\textbf{Case 3(iii):} 
If $\kappa_j<1-\mu+1/p$ for all $j\in \calJ$, then we employ the spaces related to $P_0$
and $P_1$ to derive
\begin{equation*}
 X_{\rho,\mu} = B_{p,p}^{k_{j_1}(\kappa_0 + \mu-1/p)/(1+\kappa_0 - \kappa_{j_1})}(\Gamma;F). 
\end{equation*}
For later purpose we finally note that in Case 3 the embedding
\begin{equation}\label{s-emb}
\E_{\rho,\mu} \hookrightarrow W^1_{p,\mu}\big (J;W^s_p(\Gamma;F)\big), \qquad 
  s:= \frac{k_{j_1}\kappa_0}{1+\kappa_0-\kappa_{j_1}},
\end{equation}
follows from  Proposition~3.2 of \cite{MS11} by interpolating of the spaces 
corresponding to $P_0$ and $P_1$.

We next consider the assumptions on the coefficients of the operators. For a Banach space $X$ of class $\calH\calT$, 
we write 
$$\F_{j,\mu}\big(J\times\Gamma; X\big) := W_{p,\mu}^{\kappa_j}\big(J; L_p(\Gamma;X)\big)
  \cap L_{p,\mu}\big(J; W_p^{2m\kappa_j}(\Gamma;X)\big), \qquad j\in \{0,...,m\}.$$
The following assumptions shall guarantee that each summand of the operators in 
(\ref{statmrI}) and (\ref{dynmaxregI}) maps continuously between the relevant  spaces 
described above. In view of the mapping properties of the traces and the derivatives, 
the multiplication with the coefficients has to be a
continuous map on  $\F_{j,\mu}\big(J\times\Gamma; X\big)$ for the relevant $X$.
Moreover, to perform localization 
and  perturbation, we require that the top order coefficients of the operators are bounded and
uniformly continuous.  

\begin{itemize}
 \item[\textbf{(SD)}] For $|\alpha|<2m$ we have either $2m(\mu-1/p) > 2m-1 + n/p$ and
 $a_\alpha \in \E_{0,\mu}(J\times\Omega;\calB(E)),$ or
$a_\alpha\in L_\infty\big (J\times \Omega; \calB(E)\big)$. For $|\alpha|=2m$ it holds 
$a_\alpha\in BU\!C(\ol{J}\times \ol{\Omega};\calB(E)).$ If $\Omega$ is unbounded then for 
$|\alpha|=2m$ in addition the limits $a_\alpha(t,\infty) := \lim_{|x|\ra \infty} a_\alpha(t,x)$ 
exist uniformly in $t\in \ol{J}$. 

\item[\textbf{(SB)}] Let $\calE_0 = \calB(E,F)$ and $\calE = \calB(E)$. For  $j=0,...,m$  and  
$|\beta|\leq m_j$ it holds either $b_{j\beta}\in C^{\tau_j, 2m\tau_j}\big (\ol{J}\times \Gamma; 
\calE\big )$ with some $\tau_j > \kappa_j$, or  $b_{j\beta}\in \F_{j,\mu}\big(J\times \Gamma; 
\calE\big)$ and  $\kappa_j > 1-\mu+1/p + \frac{n-1}{2mp}$.

 \item[\textbf{(SC)}] Let $\calF_0 = \calB(F)$ and $\calF= \calB(F,E)$ and 
 let $\g:V\ra \Gamma$ be any coordinates for $\Gamma$. For $j=0,...,m$  and $|\gamma|\leq k_j$ it 
 holds either $c_{j\gamma}^{\g} \in C^{\tau_j, 2m\tau_j}\big (\ol{J}\times \Gamma; \calE\big)$ 
 with some $\tau_j > \kappa_j$, or  $c_{j\gamma}^{\g}\in \F_{j,\mu}\big(J\times V;\calF\big )$
 and $\kappa_j > 1-\mu+1/p + \frac{n-1}{2mp}$.
\end{itemize}

We discuss these assumptions. First, in (SD) one can relax the boundedness assumptions
for $|\alpha|<2m$, see e.g.\ \cite{Mey10}. The fact that 
$$ C^{\tau_j, 2m\tau_j}\big (\ol{J}\times \Gamma; \calB(X)\big) 
  \cdot \F_{j,\mu}\big(J\times \Gamma;X\big) \hra \F_{j,\mu}\big(J\times \Gamma;X\big)$$
for $\tau_j > \kappa_j$ can be seen using the intrinsic norm for $W_{p,\mu}^{\kappa_j}$ and 
$W_p^{2m\kappa_j}$ given in equation~(2.6) in \cite{MS11} and Section~1 of \cite{Ama97}, 
respectively. If $\kappa_j > 1-\mu+1/p + \frac{n-1}{2mp}$, then Theorem~4.2 of \cite{MS11} and 
Sobolev's embeddings show that
$$\F_{j,\mu}\big (J\times \Gamma;X\big ) \hra BU\!C\big( \ol{J}\times \Gamma; X\big).$$
Using this fact and again the intrinsic norms of $W_{p,\mu}^{\kappa_j}$ and 
$W_p^{2m\kappa_j}$, we then derive
$$\F_{j,\mu}\big (J\times \Gamma;\calB(X)\big) \cdot \F_{j,\mu}\big (J\times \Gamma;X\big) 
   \hra \F_{j,\mu}\big (J\times \Gamma;X\big).$$
The assumption $\kappa_j > 1-\mu+1/p + \frac{n-1}{2mp}$ is only valid if $p$ and $\mu>1/p$ 
are sufficiently large. In fact, the assumption is equivalent to $p(2m\mu-m_j)>n+2m$.
The conditions in (SB) and (SC) are not optimal. For all 
$p\in (1,\infty)$, one can determine weaker Sobolev regularities for the coefficients than the ones 
given here which still meet  the requirements described above, see \cite{DHP07}, 
\cite{DPZ08} and Section~1.3.4 of \cite{Mey10}. On the other hand, (SB) and (SC) are already 
sufficient for the applications to quasilinear problems, see e.g.\  \cite{Mey10} and \cite{Mey11}.

We next state the structural assumptions on the operators, which are the same as in 
\cite{DHP07} and \cite{DPZ08}. In the sequel, the subscript $\sharp$ denotes the principle part of 
a differential operator, with an important exception for the $\calC_j$ where we put
$$\calC_{j\sharp} := 0 \qquad \text{if }\, j\notin \calJ.$$
We thus consider only the principle parts of the operators $\calC_j$ corresponding to a point on the leading 
part of the Newton polygon for $\E_{\rho,\mu}$. First, we
assume that $\calA$ is normally elliptic:

\begin{itemize}
\item[\textbf{(E)}] For all $t\in \ol{J}$, $x\in \ol{\Omega}$ and $|\xi|=1$, it holds 
$\sigma\big (\calA_\sharp(t,x,\xi)\big)\subset \C_+:=\{\rea z >0\}.$ If $\Omega$ is unbounded, 
then it holds in addition $\sigma\big(\calA_\sharp(t,\infty,\xi)\big)\subset \C_+$ for all 
$t\in \ol{J}$ and $|\xi|=1$.
\end{itemize}

We further need conditions of Lopatinskii-Shapiro type. In their formulation, local 
coordinates $\g$ for the boundary $\Gamma$ are called \textsl{associated} to $x\in \Gamma$ 
if the corresponding chart $(U,\vphi)$ satisfies 
$$\vphi(x) = 0, \qquad\vphi'(x) = \calO_{\nu(x)}, \qquad \vphi(U\cap \Omega) \subset \R_+^n, 
    \qquad \vphi(U\cap \Gamma) \subset \R^{n-1}\times \{0\},$$
where $\calO_{\nu(x)}$ is a fixed orthogonal matrix that rotates the outer normal $\nu(x)$ 
of $\Gamma$ at $x$ to $(0,...,0,-1)\in \R^n$. It is easy to see that such a chart $(U,\vphi)$ 
always exists. For coordinates $\g$ associated to $x\in \Gamma$, we define the rotated 
operators $\calA^\nu$ and $\calB_j^\nu$ by  
$$\calA^{\nu}(t,x,D): = \calA \big (t,x,\calO_{\nu(x)}^TD\big ), \qquad \calB_j^\nu(t,x,D) 
  :=  \calB_j\big (t,x, \calO_{\nu(x)}^TD\big ),\quad j=0,...,m.$$
Moreover, we introduce the local operators $\calC_j^{g}$ with respect to $\g$ by setting
$$ \calC_j^{\g}(t,x,D_{n-1}) := \sum_{|\gamma|\leq k_j} c_{j\gamma}^{\g}(t,\g^{-1}(x)) 
          D_{n-1}^\gamma, \qquad j=0,...,m,$$
where $c_{j\gamma}^{\g}$ are the local coefficients from the definition of $\calC_j$. 
We continue with the second structural assumption concerning (\ref{statmrI}).
\begin{itemize}
 \item[\textbf{(LS${}_{\text{stat}}$)}] For each fixed $t\in \ol{J}$ and $x\in \Gamma$,
  for each $\lambda\in \ol{\C_+}$ and $\xi'\in \R^{n-1}$  with $|\lambda| + |\xi'| \neq 0$ 
  and each $h\in E^m$ the ordinary initial value problem 
\begin{align*}
\lambda v(y) + \calA_\sharp^{\nu}(t, \xi', D_y)v(y) & =  0,  \qquad y>0, \\
\calB_{j\sharp}^{\nu}(t, \xi', D_y)v|_{y=0} & =  h_j,    \qquad  j=1,...,m, 
\end{align*}
has a unique solution $v\in C_0([0,\infty);E)$.
\end{itemize}
Here the space $C_0([0,\infty);E)$ consists of the continuous $E$-valued functions on 
$[0,\infty)$ vanishing at $\infty$. For the problem (\ref{dynmaxregI}) with relaxation type 
boundary conditions we need two assumptions of Lopatinskii-Shapiro type in the Cases 2 and 3.
First, in each case we require a natural  analogue of (LS${}_{\text{stat}}$).
\begin{itemize}
\item[\textbf{(LS${}_{\text{rel}}$)}] For each fixed $x\in \Gamma$, choose coordinates $\g$ 
associated to $x$. Then for every $t\in \ol{J}$,  $\lambda\in \ol{\C_+}$ and  
$\xi'\in \R^{n-1}$ with $|\lambda|+|\xi'| \neq 0$,  $h_0\in F$ and 
$h_j\in E$, $j=1,...,m$, the ordinary initial value problem 
\bea
\big (\lambda + \calA_\sharp^{\nu}(t,x, \xi', D_y)\big )v(y) & = & 0, \qquad  y>0,\\
\calB_{0\sharp}^{\nu}(t,x, \xi', D_y)v|_{y=0}+ \big ( \lambda + \calC_{0\sharp}^{\g}(t,x,\xi')\big )\sigma  & = & h_0, \\
 \calB_{j\sharp}^{\nu}(t,x, \xi', D_y)v|_{y=0}+  \calC_{j\sharp}^{\g}(t,x,\xi')\sigma  & = & h_j,\qquad j=1,...,m,
\eea  
has a unique  solution $(v,\sigma)\in C_0([0,\infty);E)\times F.$ 
\end{itemize}
In the Cases 2 and 3, the following additional `asymptotic'  conditions are required, respectively.
\begin{itemize}
\item[$\boldsymbol{(\textbf{LS}_\infty^-)}$] Let $l<2m$. For each fixed $x\in \Gamma$, choose 
coordinates $\g$ associated to $x$. Then for every $t\in \ol{J}$,  $h_0\in F$,  $h_j\in E$, 
$j=1,...,m$, and each $\lambda\in \ol{\C_+}$, $\xi'\in \R^{n-1}$ with $|\lambda|+|\xi'| \neq 0$,   
the ordinary initial value problem 
\bea
\big (\lambda + \calA_\sharp^\nu(t,x, \xi', D_y)\big )v(y) & = & 0, \qquad  y>0,\\
 \calB_{j\sharp}^\nu(t,x, \xi', D_y)v|_{y=0} & = & h_j,\qquad j=1,...,m,
\eea  
and for all $\lambda\in \ol{\C_+}$ and $|\xi'| = 1$ the problem
\bea
 \calA_\sharp^\nu(t,x, \xi', D_y)v(y) & = & 0, \qquad  y>0,\\
\calB_{0\sharp}^\nu(t,x, \xi', D_y)v|_{y=0} + \big ( \lambda + \calC_{0\sharp}^{\g}(t,x,\xi')\big )\sigma  & = & h_0, \\
 \calB_{j\sharp}^\nu(t,x, \xi', D_y)v|_{y=0} +  \calC_{j\sharp}^{\g}(t,x,\xi')\sigma  & = & h_j,\qquad j=1,...,m,
\eea  
has a unique solution $(v,\sigma)\in C_0([0,\infty);E)\times F$, respectively.
\item[$\boldsymbol{(\textbf{LS}_\infty^+)}$] Let $l>2m$.  For each fixed $x\in \Gamma$, choose 
coordinates $\g$ associated to $x$. Then for each $t\in \ol{J}$,  $h_0\in F$,  $h_j\in E$, 
$j=1,...,m$, and each $\lambda\in \ol{\C_+}$ and $\xi\in \R^{n-1}\bs\{0\}$, 
the ordinary initial value problem 
\bea
\big (\lambda + \calA_\sharp^\nu(t,x, \xi', D_y)\big )v(y) & = & 0, \qquad  y>0,\\
 \calB_{j\sharp}^\nu(t,x, \xi', D_y)v|_{y=0} + \delta_{j,\, \calJ_{2q_{\text{max}}+1}}\calC_{j\sharp}^{\g}(t,x,\xi')\sigma & = & h_j,\qquad j=0,...,m,
\eea  
and further for all $\lambda\in \ol{\C_+}\bs\{0\}$, $|\xi'| = 1$  and $q=1,..., 2q_{\text{max}}$, the problem
\bea
\big (\lambda + \calA_\sharp^\nu(t,x, 0, D_y)\big )v(y) & = & 0, \qquad  y>0,\\
\calB_{0\sharp}^\nu(t,x, 0, D_y)v|_{y=0} + \delta_{-1,\, \calJ_q}\lambda\sigma + \delta_{0,\,\calJ_q} \calC_{0\sharp}^{\g}(t,x,\xi')\sigma  & = & h_0, \\
 \calB_{j\sharp}^\nu(t,x, 0, D_y)v|_{y=0} +  \delta_{j,\,\calJ_q} \calC_{j\sharp}^{\g}(t,x,\xi')\sigma  & = & h_j,\qquad j=1,...,m,
\eea 
has a unique solution $(v,\sigma)\in C_0([0,\infty);E)\times F$, respectively.
Here  $\delta_{j,\,\calJ_q} = 1$ if $j\in \calJ_q$ and $\delta_{j,\,\calJ_q} = 0$ otherwise.
\end{itemize}
In \cite{DHP07} and \cite{DPZ08} it is shown that these conditions are necessary for 
 maximal $L_p$-regularity of (\ref{dynmaxregI}) on finite intervals. In Section~3 of 
 \cite{DPZ08} they are verified for a variety of concrete problems from the applications,
 see also \cite{Mey10} and \cite{Mey11}. If $E$ and $F$ are finite dimensional, it suffices 
 to consider $h_0=h_j=0$ in  the above conditions.

We can now state our maximal $L_{p,\mu}$-regularity results. We start with the one for 
(\ref{statmrI}).

\begin{thm}\label{statmrthm}\textsl{Let $E$ be a Banach space of class $\HT$, $p\in (1,\infty)$ and 
$\mu\in (1/p,1]$. Let $J=(0,T)$ be a finite interval, and let $\Omega\subset \R^n$ be a domain with 
compact smooth boundary $\Gamma=\partial\Omega$. Assume that \emph{(E), (LS${}_{\text{stat}}$), (SD)}
 and \emph{(SB)} hold true and that $\kappa_j \neq 1-\mu+1/p$ for $j=1,...,m$. Then the problem
\begin{alignat}{3}
\partial_t u + \calA(t,x,D) u & =    f(t,x), & \qquad & x\in  \Omega, & \qquad t\in J,&   \nonumber\\
\calB_j(t,x,D) u & =  g_j(t,x), && x\in \Gamma,  & t\in J,&  \qquad j=1,...,m, \nonumber\\
u(0,x) & =  u_0(x), &&  x\in \Omega, &&  \nonumber
\end{alignat}
has a unique solution $u = \calL{}_{\text{\emph{stat}}}(f,\wt{g},u_0)\in \E_{u,\mu}$ if and only if 
$(f,\wt{g},u_0)\in  \calD_{\emph{\text{stat}}}$, where 
\begin{align*}
\calD_{\emph{\text{stat}}} := \big \{(f,\wt{g},  u_0)\in  \E_{0,\mu}\times \wt{\F}_{\mu}&
         \times X_{u,\mu}\;:\;\text{ for $j=1,...,m$ it holds}\\
&\,\calB_j(0,\cdot,D) u_0 = g_j(0,\cdot)  \text{ on } \Gamma \text{ if } \kappa_j >1-\mu-1/p\big\}.
\end{align*}
The corresponding solution operator $\calL{}_{\text{\emph{stat}}}:\calD_{\emph{\text{stat}}}\ra 
\E_{u,\mu}$ is continuous. If $\calL{}_{\text{\emph{stat}}}$ is restricted to 
\begin{align*}
\calD_{\emph{\text{stat}}}^0:= \big \{ (f,\wt{g},u_0) \in \calD_{\emph{\text{stat}}}\;:\; 
   g_j|_{t=0} =0 \text{ if } \kappa_j >1-\mu-1/p,\;\; j=1,...,m\big \},
\end{align*} 
for any given $T_0>0$ the operator norm of the restriction is uniformly bounded 
for  $T \in(0, T_0]$.}
\end{thm}

In the situation of the theorem, it is clear that if the coefficients
$$
(-\ii)^{|\alpha|}a_\alpha, \quad |\alpha|\leq 2m, \qquad (-\ii)^{|\beta|}b_{j\beta}, \quad |\beta|\leq m_j, \quad j=1,...,m,
$$
and the data are real-valued, then also the solution $u$ is real-valued. We next state
the maximal regularity result for (\ref{dynmaxregI}).

\begin{thm}\label{dynmrthm}\textsl{Let $E$ and $F$ be Banach spaces of class $\HT$, 
$p\in (1,\infty)$ and $\mu\in (1/p,1]$. Let $J=(0,T)$ be a finite interval,  and let 
$\Omega\subset \R^n$ be a  domain with compact smooth boundary $\Gamma=\partial\Omega$. Assume 
that \emph{(E)}, \emph{(LS${}_{\text{rel}}$)}, \emph{(SD)}, \emph{(SB)} and \emph{(SC)} are valid 
and that, in addition, if $l<2m$ condition \emph{(LS${}_\infty^-$)} holds and if $l>2m$ 
condition \emph{(LS${}_\infty^+$)} holds. Assume further that $\kappa_j \neq 1-\mu+1/p$ for 
$j=0,...,m$. Then the problem
\begin{alignat}{3}
\partial_t u + \calA(t,x,D) u & =    f(t,x), & \qquad & x\in  \Omega, & \quad &  t\in J,  \nonumber\\
\partial_t \rho + \calB_0(t,x,D) u + \calC_0(t,x,D_\Gamma) \rho & =  g_0(t,x), &&x\in  \Gamma,  &&  t\in J,\quad   \nonumber\\
\calB_j(t,x,D) u + \calC_j (t,x,D_\Gamma)\rho & =  g_j(t,x), &&x\in  \Gamma,  &&  t\in J,\qquad j=1,...,m,   \nonumber\\
u(0,x) & =  u_0(x), &&  x\in  \Omega, &&  \nonumber \\
\rho(0,x) & =  \rho_0(x), &&  x\in  \Gamma, &&  \nonumber
\end{alignat}
has a unique solution $(u,\rho)\in \E_{u,\mu}\times \E_{\rho,\mu}$ if and only if 
$(f,g,u_0,\rho_0)\in \calD_{\emph{\text{rel}}}$, where 
\begin{align*}
 \calD_{\emph{\text{rel}}}:= \big \{ (f,g,&u_0,\rho_0)\in \E_{0,\mu}\times \F_{\mu} \times X_{u,\mu} \times X_{\rho,\mu}\;:\;\;\; \text{for $j=1,...,m$ it holds}\\
 &\,\calB_j(0,\cdot,D)u_0 + \calC_j(0,\cdot ,D_\Gamma)\rho_0 = g_j(0,\cdot ) \; \text{ on }\,  \Gamma \;\;\text{ if } \kappa_j>1-\mu+1/p;\\
\;\;\qquad &\,\quad  g_0(0,\cdot) - \calB_0(0,\cdot,D)u_0 - \calC_0(0,\cdot,D_\Gamma) \rho_0 \in X_{\partial_t \rho,\mu} \qquad \text{if }\; \kappa_0>1-\mu+1/p \big \}.
\end{align*}
The corresponding solution operator $\calL_{\emph{\text{rel}}}: \calD_{\emph{\text{rel}}} \ra 
\E_{u,\mu}\times \E_{\rho,\mu}$ is continuous. If $\calL_{\emph{\text{rel}}}$ is restricted to 
\begin{align*}
  \calD_{\emph{\text{rel}}}^0 := \big \{ (f,g,u_0,\rho_0) \in  \calD_{\emph{\text{rel}}} \;:\; g_j|_{t=0} =0 \text{ if } \kappa_j >1-\mu-1/p,\;\; j=0,...,m\big \},
\end{align*} 
for any given $T_0>0$ the operator norm of the restriction is uniformly bounded for $T\in(0, T_0]$. }
\end{thm}

A similiar statement as above holds for real-valued solutions.  In the theorems, the spaces 
$\calD_{\text{stat}}$ and $\calD_{\text{rel}}$  
are considered as closed subspaces of  $\E_{0,\mu}\times \wt{\F}_{\mu}\times X_{u,\mu}$
and $ \E_{0,\mu}\times \F_{\mu} \times X_{u,\mu} \times X_{\rho,\mu}$, respectively.
They contain the compatibility conditions of the boundary inhomogeneities and the initial values at 
$t=0$, which are necessary for the solvability of (\ref{statmrI}) and  (\ref{dynmaxregI}), 
respectively. 

One needs the spaces $\calD_{\text{stat}}^0$ and $\calD_{\text{rel}}^0$ 
with vanishing initial values since the resulting uniform estimates are crucial for 
fixed point arguments arising in the context of  quasilinear problems. They are considered as closed subspaces of $\E_{0,\mu}\times {}_0 \wt{\F}_{\mu}\times X_{u,\mu}$
and $ \E_{0,\mu}\times {}_0\F_{\mu} \times X_{u,\mu} \times X_{\rho,\mu}$, respectively, where  ${}_0\wt{\F}_{\mu}$  and ${}_0 \F_\mu$ are defined as follows. For a Banach space $X$ of class $\calH\calT$ and $s= [s] + s_*$ with $[s]\in \N_0$, $s_*\in [0,1)$ we set
$${}_0W_{p,\mu}^s(J; X) := \big ({}_0W_{p,\mu}^{[s]}(J;X), {}_0W_{p,\mu}^{[s]+1}(J;X)\big)_{s_*,p},$$
where ${}_0W_{p,\mu}^k(J;X) := \big \{u\in W_{p,\mu}^k(J;X):u(0), ..., u^{(k-1)}(0) = 0 \big \}$ is considered as a closed subspace of  $W_{p,\mu}^k(J;X)$, and then 
$${}_0 \F_{j,\mu} := {}_0W_{p,\mu}^{\kappa_j}\big (J; L_p(\Gamma;E)\big) \cap 
          L_{p,\mu}\big(J;W_p^{2m\kappa_j}(\Gamma;E)\big),\qquad j=1,...,m.$$ 
Analogously we define the spaces ${}_0 \F_{0,\mu}$, ${}_0 \wt{\F}_{\mu}$, ${}_0 \F_{\mu}$, ${}_0\E_{u,\mu}$ and 
${}_0\E_{\rho,\mu}$.  It is shown in Proposition~2.10 of \cite{MS11} that ${}_0W_{p,\mu}^s = W_{p,\mu}^s$ if
$0<s<1-\mu+1/p$ and 
$${}_0W_{p,\mu}^s = \big\{ u\in W_{p,\mu}^s\;:\; u^{(l)}(0) = 0,\;\; 
        l\in \{0,...,k\}\big\}$$
if $k+1-\mu+1/p < s < k + 2-\mu+1/p$ for $k\in\N_0$, with equivalent norms, respectively. In other words, the trace at $t=0$ of a 
derivative of $u\in {}_0W_{p,\mu}^s$ vanishes if it exists.

Compared to the unweighted approach, the maximal $L_{p,\mu}$-regularity approach has the following advantages, where we restrict to the setting of (\ref{statmrI}). Analogous statements are valid for (\ref{dynmaxregI}).
\begin{itemize}
 \item \textsl{Flexible initial regularity:} We obtain solutions for initial values in $B_{p,p}^s(\Omega;E)$, where $s\in (0, 2m(1-1/p)]$.
 \item \textsl{Inherent smoothing effect:}  Away from the initial time, $\tau\in (0,T)$, the solutions belong to   
$$ W_{p}^1\big (\tau,T; L_p(\Omega;E)\big )\cap L_{p}\big (\tau,T; W_p^{2m}(\Omega;E)\big )\hra C\big (\ol{J}; B_{p,p}^{2m(1-1/p)}(\Omega;E)\big ).$$ 
\item \textsl{Control solutions in a strong norm at a later time  by a weaker norm at an earlier time and the data}: For $s= 2m(\mu-1/p) \in (0, 2m(1-1/p)]$ it holds $$|u(T)|_{B_{pp}^{2m(1-1/p)}(\Omega;E)} \leq C(T)(|f|_{\E_{0,\mu}} + |\wt{g}|_{\wt{\F}_{\mu}} + |u_0|_{B_{p,p}^s(\Omega;E)}).$$
\item \textsl{Avoid compatibility conditions:} Given $p\in (1,\infty)$, if $\mu$ is sufficiently 
close to $1/p$, we have $\kappa_j < 1-\mu+1/p$ for all $j$ and thus obtain  a unique solution 
$u\in \E_{u,\mu}$ for arbitrary data in $\E_{0,\mu}\times \wt{\F}_{\mu} \times X_{u,\mu}$.
\end{itemize}

The rest of the paper is concerned with the proofs of the Theorems \ref{statmrthm} and 
\ref{dynmrthm}.

\section{The model problems}
We first consider the full-space case $\Omega = \R^n$ without boundary conditions and assume that 
the coefficients of the differential operator 
$$\calA(D) =\sum_{|\alpha|=2m}a_\alpha D^\alpha$$
are constant, i.e., $a_\alpha \in \calB(E)$ are independent of $(t,x)$. Observe that there are no 
lower order terms and that $\calA$ is homogeneous of degree $2m$. We have the following maximal 
$L_{p,\mu}$-regularity result for $\calA$ on the half-line.

\begin{lem}\label{sec:HSP2easy} \textsl{Let $E$ be a Banach space of class $\HT$, $p\in (1,\infty)$, 
$\mu\in (1/p,1]$, and assume that the constant coefficient operator $\calA$ satisfies \emph{(E)}. 
Then there is a unique solution $u = \calS_F(f,u_0)\in \E_{u,\mu}(\R_+\times \R^n)$ of 
\begin{alignat}{3}
u + \partial_tu + \calA(D) u & =  f(t,x), & \qquad & x\in \R^n, \qquad t>0, \nonumber\\
u(0,x) & =  u_0(x), && x\in \R^n,   \label{eq:HSP1}
\end{alignat}
if and only if 
$$f\in \E_{0,\mu}(\R_+\times \R^n) \quad \text{ and } \quad  u_0\in X_{u,\mu}(\R^n).$$}
\end{lem}
\bprf Lemma~4.2 of \cite{DHP07} shows that the realization of the operator $1+\calA$ on $L_p(\R^n;E)$
 with domain $D(1+\calA)= W_p^{2m}(\R^n;E)$ admits maximal $L_p$-regularity on the half-line. Since 
$$X_{u,\mu}(\R^n) = B_{p,p}^{2m(\mu-1/p)}(\R^n;E) 
     =  \big (L_p(\R^n;E), W_p^{2m}(\R^n;E)\big )_{\mu-1/p,p},$$ 
the assertion follows from Theorem~3.2 of \cite{PS04}.\eprf

The model problems for (\ref{statmrI}) and (\ref{dynmaxregI}) on the half-space involve 
 boundary conditions and thus require a much greater effort. To construct the solution,
one uses an operator-valued Fourier multiplier theorem in $L_{p,\mu}$. For Banach spaces $X$, $Y$ 
and a symbol $m\in L_{1,\text{loc}}\big (\R; \calB(X,Y)\big )$ one 
introduces an operator $T_m$ by setting 
$$T_m f := \calF^{-1} m \calF f, \qquad f\in \calF^{-1} C_c^\infty(\R;X),$$ 
where $\calF$ denotes the Fourier transform on the real line. We can restrict $T_m$ 
to functions on $\R_+$. Observe that 
$T_m$ is densely defined on $L_{p,\mu}(\R_+;X)$. We also use the analogous definition
on the space $L_{p}(\R^n;X)$. The next result is due to Girardi $\&$ Weis \cite{GW11}. 

\begin{thm}\label{sec:mumult}\textsl{Let $p\in (1,\infty)$, $\mu\in (1/p,1]$, and let $X, Y$ be 
Banach spaces of class $\calH\calT$. Assume that  $m\in C^1\big (\R\bs \{0\}; \calB(X,Y)\big)$ 
satisfies $\calR\big(\{m(\lambda), \lambda m'(\lambda) : \lambda\neq 0\}\big) <\infty.$ 
Then $T_m \in \calB\big (L_{p,\mu}(\R_+;X), L_{p,\mu}(\R_+;Y)\big )$.}
\end{thm}
Here, the $\calR$-bound of a family $\calT\subset \calB(X,Y)$ is denoted by $\calR(\calT)$. 
For a definition and properties of $\calR$-boundedness we refer to \cite{DHP03} or \cite{KW04}.
Under more restrictive assumptions on the symbol $m$ we can give a short proof a multiplier theorem 
in $L_{p,\mu}$, employing a result of Kr\'ee \cite{Kr66} (which is also used in the proof 
in \cite{GW11}).

\begin{prop}\label{sec:mumultprop}\textsl{In addition to the assumptions of 
Theorem \ref{sec:mumult}, suppose that  $m$ satisfies 
$$m\in C^2\big (\R\bs \{0\}; \calB(X,Y)\big), \qquad |m''(\lambda)|_{\calB(X,Y)} 
 \lesssim |\lambda|^{-2} \;\text{ for }\;\lambda \neq 0.$$ 
Then $T_m \in \calB\big (L_{p,\mu}(\R_+;X), L_{p,\mu}(\R_+;Y)\big )$.}
\end{prop}
\bprf The operator-valued multiplier theorem for the unweighted case $\mu=1$ shows that $T_m$ 
extends to a bounded operator from $L_p(\R_+;X)$ to $L_p(\R_+;Y)$; see Theorem~3.4 of \cite{Wei01}. 
Moreover, following the lines of the proof of Lemma~VI.4.4.2 of \cite{St93}, the assumptions on $m$ 
imply that $T_m$ may be represented as a convolution  operator with a kernel 
$k\in C\big (\R\bs\{0\}; \calB(X,Y)\big )$ satisfying $|k(t)|_{\calB(X,Y)} \lesssim |t|^{-1}$. 
It now follows from  Th\'eor\`eme~2 of \cite{Kr66} that $T_m$ is also bounded from 
$L_{p,\mu}(\R_+;X)$ to $L_{p,\mu}(\R_+;Y)$, for all $\mu\in (1/p,1]$.\eprf

We next treat the half-space model problem corresponding to (\ref{dynmaxregI}),
where we proceed similarly as in Section~4 of \cite{DPZ08}.
On  $\Omega = \R_+^n$ with boundary $\Gamma = \R^{n-1}$ we consider the homogeneous 
differential operator 
$$\calA(D) =\sum_{|\alpha|=2m}a_\alpha D^\alpha$$ 
and  the homogeneous boundary operators 
$$\calB_j(D) = \sum_{|\beta|= m_j} b_{j\beta} \tr_{\R_+^n} D^\beta, \qquad  \calC_j(D_{n-1}) = \sum_{|\gamma|= k_j} c_{j\gamma}  D_{n-1}^\gamma \qquad j=0,...,m.$$ 
The coefficients of the operators  
$$a_\alpha, b_{j\beta}\in \calB(E), \quad c_{j\gamma}\in \calB(F,E), \quad j=1,...,m,
  \quad b_{0\beta}\in \calB(E,F) ,\quad c_{0\gamma}\in \calB(F)$$ 
  are assumed to be independent of $t$ and $x$.
If nothing else is indicated, now all spaces have to be understood over $\R_+\times \R_+^n$ 
and over $\R_+\times \R^{n-1}$.

\begin{lem}\label{lg10}\textsl{Let $E$ and $F$ be Banach spaces of class $\HT$, $p\in (1,\infty)$, 
and $\mu\in (1/p,1]$. We assume that \emph{(E)} and \emph{(LS${}_{\text{rel}}$)} are valid and
that condition \emph{(LS${}_\infty^-$)} holds if $l<2m$ and condition \emph{(LS${}_\infty^+$)}
holds if $l>2m$. Let $(f,g,u_0,\rho_0)\in \calD_{\emph{\text{rel}}}$. Then there is a unique 
solution  $(u,\rho)\in \E_{u,\mu}\times \E_{\rho,\mu}$ of
\begin{alignat}{3}
u + \partial_t u + \calA(D)u & =    f(t,x), & \qquad & x\in \R_+^n, & \quad &  t>0,  \nonumber\\
\rho + \partial_t \rho + \calB_0(D) u + \calC_0(D_{n-1}) \rho & =  g_0(t,x), &&x\in \R^{n-1},  &&  t>0,\quad   \nonumber\\
\calB_j(D) u + \calC_j(D_{n-1}) \rho & =  g_j(t,x), &&x\in \R^{n-1},  &&  t>0,\quad j=1,...,m,  \label{dynmaxregHS}\\
u(0,x) & =  u_0(x), &&  x\in \R_+^n, &&  \nonumber \\
\rho(0,x) & =  \rho_0(x), &&  x\in \R^{n-1}. &&  \nonumber
\end{alignat}}
\end{lem}
\bprf \textbf{(I)} We first show uniqueness for  (\ref{dynmaxregHS}). We use the space
$Z:=L_p(\R_+^n;E) \times W_p^s(\R^{n-1};F)$ with
 $s= 2m\kappa_0$ in the Cases 1 and 2 as well as $s= k_{j_1}\kappa_0/(1+\kappa_0-\kappa_{j_1})$ 
in Case 3. On $Z$, we introduce the operator $A$ defined by
$$A(u,\rho):= \big ((1+ \calA)u,  \calB_0u + (1+\calC_0)\rho\big), \qquad (u,\rho)\in D(A),$$ 
with domain
\begin{align*}
 D(A) := \big \{ (u,\rho)  \in  W_p^{2m}&(\R_+^n;E)\times W_p^{l+2m\kappa_0}(\R^{n-1};F)\;:\;\\
&\,  \calB_j u + (1+\calC_j) \rho = 0,\;\;j=1,...,m;\;\;\; \calB_0 u + \calC_0\rho \in W_p^s(\R^{n-1};F)\big \}.
\end{align*}
By (the proof of) Theorem~2.2 of \cite{DPZ08}, $A$ generates an analytic $C_0$-semigroup 
on $Z$. Due to \eqref{E-rho} and \eqref{s-emb}, the space  $\E_{u,\mu} \times \E_{\rho,\mu}$
embeds into 
$$\mathbb{G}:=\E_{u,\mu}(\R_+) \times \big ( W_{p,\mu}^1(\R_+; W_p^s(\R^{n-1};F)) 
 \cap L_{p,\mu}(\R_+; W_p^{l+2m\kappa_0}(\R^{n-1};F))\big ).$$ 
Let $u\in\mathbb{G}$ be a solution of  \eqref{dynmaxregHS} with $u_0=0$, $\rho_0=0$,
$f=0$ and $g_0=\cdots=g_m=0$.  Since $L_{p,\mu}(J;Z)\hra L_1(J;Z)$, it follows that $u$ 
is a mild solution of the inhomogeneous evolution equation for $A$ on $Z$ with trivial data, 
and thus $u=0$.

\textbf{(II)} The rest of the proof is concerned with the existence of solutions of 
(\ref{dynmaxregHS}). We write $x= (x',y)\in \R_+^n$ with $x'\in \R^{n-1}$ and $y>0$, as well
as $\calF_{x'}$ and $\calF_t$ for the partial Fourier transform with respect to $x'$
and $t\in\R$, with covariable  $\xi'\in \R^{n-1}$ and $\theta\in \R$, respectively. 
In order to apply  $\calF_t$, we extend a function with compact support in $\R_+$
by 0 to $\R$. In the same way as in Section~4.1 of \cite{DPZ08} one can see that 
it sufficies to consider the case
$$f=0, \qquad g= (g_0,..., g_m)\in {}_0\F_{\mu}, \qquad u_0 = 0, \qquad \rho_0 = 0.$$ 
(See Lemma~3.2.2 and Proposition~3.2.3 of \cite{Mey10}.)
Moreover we first assume that
$$g \in \calD:=C_c^\infty(\R_+\times \R^{n-1};F\times E^m). $$
It can be seen as in Lemma~1.3.14 of \cite{Mey10} that $\calD$ is dense in ${}_0\F_{\mu}$.
For such data the problem (\ref{dynmaxregHS}) was solved in the proof of Theorem~2.1 of \cite{DPZ08}.
In the following we estimate the norm of the solution $(u,\rho)$ in  the weighted
solution space ${}_0\E_{u,\mu} \times {}_0\E_{\rho,\mu}$ by the norm of $g$ in ${}_0\F_{\mu}$. For this 
estimate, we have to derive an appropriate representation of $(u,\rho)$.
We apply  $\calF_{x'}\calF_t$  to (\ref{dynmaxregHS}) and arrive for any $\theta\in \R$ 
and $\xi'\in \R^{n-1}$ at the ordinary initial value problem
\begin{alignat}{3}
(1+\ii\theta)v + \calA(\xi', D_y) v & =    0, & \qquad & y>0,  \nonumber\\
(1+\ii\theta)\sigma +  \calB_0(\xi', D_y) v |_{y=0} + \calC_0(\xi') \sigma & = \big(\calF_{x'}\calF_tg_0\big)(\theta, \xi'),  
      && &&  \quad   \label{jesaba36}\\
\calB_j(\xi', D_y) v|_{y=0}+ \calC_j(\xi') \sigma & =  \big(\calF_{x'}\calF_tg_j\big)(\theta, \xi'), && j=1,...,m. 
   &&     \nonumber
\end{alignat}
In Section~4.3 of \cite{DPZ08} it is shown that (\ref{jesaba36}) possesses for all $\theta$ and $\xi'$ 
a unique solution $\big (v(\theta,\xi', \cdot), \sigma(\theta, \xi')\big)$ which may be represented as follows.  
We define the symbols  
$$\vartheta := (1+\ii\theta+ |\xi'|^{2m})^{1/2m}, \qquad  b := \frac{|\xi'|}{\vartheta}, 
  \qquad \zeta := \frac{\xi'}{|\xi'|}, \qquad a:= \frac{1+\ii\theta}{\vartheta^{2m}},$$ 
and the so-called boundary symbol $s(\theta,\xi')$ by
\begin{align*}
s(\theta,\xi') &:= 1+\ii\theta + |\xi'|^{l}\qquad \text{ in the Cases 1 and 2},\\
s(\theta,\xi') &:= 1+\ii\theta + \sum_{j\in \calJ} |\xi'|^{k_j} \vartheta^{m_0-m_j}\qquad 
      \text{ in  Case 3}.
\end{align*} 
Then it holds
\begin{align*}
v(\theta,\xi',y) &= \text{ first component of }\;  e^{\vartheta 
          \ii A_0(b\zeta,a)y} P_s(b\zeta,a) M_u^0(b,\zeta,\vartheta)\big(\vartheta^{-m_j}\calF_{x'}\calF_t g_j(\theta, \xi')\big)_{j=0,...,m},\\
\sigma(\theta,\xi') &=   s(\theta,\xi')^{-1}  \vartheta^{m_0}   
      M_{\rho}^0(b,\zeta,\vartheta) \big(\vartheta^{-m_j}\calF_{x'}\calF_t g_j(\theta, \xi')\big)_{j=0,...,m}.
\end{align*} 
Here we have used  holomorphic functions
$$A_0:\C^{n-1}\times \C\ra  \calB(E^{2m}), \qquad P_s:\C^{n-1}\times \C\ra\calB(E^{2m}),$$
$$M_u^0: D_{b} \times D_\zeta \times \Sigma\ra  \calB(F\times E^m, E^{2m}), 
\qquad M_\rho^0: D_{b} \times D_\zeta   \times \Sigma\ra \calB(F\times E^m, F),$$ 
where  $D_b \subset \C$ and $D_\zeta \subset \C^{n-1}\setminus\{0\}$ are bounded open sets 
satisfying
$$(\ol{B}_{1/2}(1/2))^{1/2m} \subset D_b, \qquad 
       \{\zeta\in \R^{n-1}\;:\; |\zeta| = 1\} \subset D_\zeta,$$ 
and $\Sigma\!=\!\Sigma_\phi \!=\! \big \{z\in \C\bs\{0\}: |\text{arg}z|< \phi\big \}$ 
is a sector with
$\phi\in (\frac{\pi}{4m},\pi)$.  The spectrum of $\ii A_0(b\zeta,a)$ has a gap at the imaginary 
axis,  and $P_s(b\zeta,a)$ is the spectral projection corresponding to the stable part of the 
spectrum.   The maps $M_u^0$ and $M_\rho^0$ have the crucial property that
\beq\label{lg7}
\big \{|\xi'|^{|\alpha'|} D_{\xi'}^{\alpha'} M_u^0(\wt{b},\xi'|\xi'|^{-1},\wt{\vartheta})
  \,:\, \alpha'\in \{0,1\}^{n-1},\;\;\xi'\neq 0, \;\; \wt{b}\in D_b, 
   \;\; \wt{\vartheta}\in \Sigma\big \}
\eeq
is an $\calR$-bounded set of operators in $\calB(F\times E^m, E^{2m})$, and that
\beq\label{kv38}
\big \{|\xi'|^{|\alpha'|} D_{\xi'}^{\alpha'} M_\rho^0(\wt{b},\xi'|\xi'|^{-1},\wt{\vartheta})\,:\, 
\alpha'\in \{0,1\}^{n-1},\;\;\xi'\neq 0, \;\; \wt{b}\in D_b, \;\; 
  \wt{\vartheta}\in \Sigma\big \}
\eeq
is an $\calR$-bounded set of operators in $\calB(F\times E^m, F)$. For the solvability and the 
representation of the solution of (\ref{jesaba36}) in \cite{DPZ08} only the condition 
(LS${}_\text{rel}$) is needed. In the Cases 2 and 3 the asymptotic 
Lopatinskii-Shapiro conditions (LS${}_\infty^-$) and (LS${}_\infty^+$) are required 
to show the $\calR$-boundedness of the sets in \eqref{lg7} and \eqref{kv38}, 
because of the unboundedness of $\vartheta$. 
In Case 1 the symbols $M_u^0$ and $M_\rho^0$ do not depend on $\vartheta$, 
so that in this case additional conditions are not necessary.

Since $\calF_{x'}\calF_tg$ belongs for $g\in \calD$ to the Schwartz class and all derivatives of 
the terms involved in the representation of the solution grow at most polynomially, we can
apply the inverse Fourier transforms and obtain that
\begin{align*}
u &= \text{ first component of }\; \calF_{t}^{-1}\calF_{x'}^{-1}  e^{\vartheta 
          \ii A_0(b\zeta,a)y} P_s(b\zeta,a) M_u^0(b,\zeta,\vartheta)  \big(\vartheta^{-m_j}\calF_{x'}\calF_{t} g_j\big)_{j=0,...,m},\\
\rho &=   \calF_{t}^{-1}\calF_{x'}^{-1} s(\theta,\xi')^{-1}  \vartheta^{m_0}   
      M_{\rho}^0(b,\zeta,\vartheta) \big(\vartheta^{-m_j} \calF_{x'}\calF_{t} g_j\big)_{j=0,...,m}
\end{align*} 
is the unique solution of (\ref{dynmaxregHS}) with $f=0$, $u_0=0$, $\rho_0=0$ and 
$g\in\calD$.

\textbf{(III)} We derive another representation of the solution 
by identifying the Fourier multipliers with operators.
For a function $h\in \calS( \R^{n-1};E^{2m})$  and fixed $(x',y)\in \R_+^n$
we calculate  
\begin{align}
 \big(  \calF_{x'}^{-1}e^{\ii \vartheta A_0y}P_sh\big )({x'}) &\,  =   \big( \calF_{x'}^{-1}e^{\ii\vartheta A_0(y+\wt{y})}P_s e^{-\wt{y}\vartheta}h\big)({x'})|_{\wt{y}=0}  \label{kd205}\\
& \,= - \int_0^\infty \partial_{\wt{y}} \big( \calF_{x'}^{-1}e^{\ii\vartheta A_0(y+\wt{y})}P_s e^{-\wt{y}\vartheta}h\big)({x'})  \D \wt{y}  \nonumber\\
& \,=   \int_0^\infty \big( \calF_{x'}^{-1} e^{\ii\vartheta A_0(y+\wt{y})}P_s \frac{1-\ii A_0}{\vartheta^{2m-1}}\vartheta^{2m} e^{-\wt{y}\vartheta}  h\big )({x'})  \D \wt{y}  \nonumber\\
& \, =  \int_0^\infty \big( \calF_{x'}^{-1} e^{\ii\vartheta A_0(y+\wt{y})}P_s \frac{1-\ii A_0}{\vartheta^{2m-1}}\big) * \big((L_{\theta} \calE_\theta\calF_{x'}^{-1}h)(\cdot,\wt{y})\big)(x') \D \wt{y}, \nonumber
\end{align}
neglecting the arguments of $A_0$ and $P_s$.
Here the operator $L_\theta$ is defined by
$$L_\theta := 1+\ii\theta + (-\Delta_{n-1})^m = \calF_{x'}^{-1}\vartheta^{2m}\calF_{x'} ,$$ 
where the last equality holds, e.g., on Schwartz functions. We observe that for a bounded
holomorphic scalar function $\varphi$ on a sector $\Sigma_\tau$ with $\tau\in(0,\pi)$
the operator $\varphi(-\Delta_{n-1})$ defined via
the $\calH^\infty$-calculus for $-\Delta_{n-1}$ on $L_p(\R^{n-1};E)$ coincides
with the Fourier multiplier $\calF_{x'}^{-1}\varphi(|\cdot|^2)\calF_{x'}$, see 
Example~10.2 of \cite{KW04}. Moreover, the $\calH^\infty$-calculus  extends the usual Dunford 
type calculus for sectorial operators, see Remark~9.9 of \cite{KW04}.  Therefore,
the extension operator $\calE_{\theta}$, which corresponds to $y\mapsto e^{- y\vartheta}$, 
is  given by
$$(\calE_{\theta}f)(x',y) := e^{-y L_{\theta}^{1/2m}} f(x'), \qquad x'\in \R^{n-1}, \qquad y >0,$$
for $f\in L_p(\R^{n-1};E)$. We  also obtain the equality
$$\calF_{x'}^{-1} \vartheta^{2m} e^{-\cdot \vartheta}h = L_{\theta} \, \calE_\theta\,\calF_{x'}^{-1}h, \qquad h\in \calS(\R^{n-1};E^m),$$
which we have used in the last line of (\ref{kd205}). For $\theta\in \R$ and $f\in L_p(\R_+^n;E^{2m})$ we thus define the operator $\calT(\theta)$ by 
\begin{align*}
(\calT(\theta)f)(x',y)  := \text{ first component of } \int_0^\infty \big( \calF_{x'}^{-1} e^{\ii\vartheta A_0(y+\wt{y})}P_s \frac{1-\ii A_0}{\vartheta^{2m-1}}\big) *f(\cdot,\wt{y}\big)(x') \D \wt{y}.
\end{align*} 
The proofs of Lemmas~4.3 and 4.4 in \cite{DHP07} show that 
$\calT\in C^1\big( \R; \calB(L_p(\R_+^n;E^{2m}), W_p^{2m}(\R_+^n;E))\big)$ and that 
\beq\label{lg6}
\big \{  D^\alpha \calT(\theta),\; \theta \frac{\partial}{\partial\theta} D^\alpha \calT(\theta)\;:\; \theta\in \R, \;\; |\alpha|\leq 2m\big \}
\eeq
is an $\calR$-bounded set of operators in $\calB\big( L_p(\R_+^n;E^{2m}), L_p(\R_+^n;E)\big)$. Further, as above one can see that $\vartheta^{-m_j} \calF_{x'} = \calF_{x'} L_\theta^{-m_j/2m}$ on Schwartz functions, for $j=0,...,m$.
This fact leads to
$$u =   \calF_{t}^{-1}\calT(\theta)  L_\theta \calE_\theta  \calF_{x'}^{-1} M_u(b,\zeta,\vartheta) \calF_{x'}\big(L_\theta^{-m_j/2m} \calF_t g_j\big)_{j=0,...,m}.$$
The Dunford type calculus for sectorial operators yields for $\theta\in \R$ and $y>0$ the representation
$$
L_{\theta} e^{- y L_{\theta}^{1/2m}} = \frac{1}{2\pi\ii} \int_{\Xi} ze^{-y z^{1/2m}} 
  (z - L_{\theta})^{-1} \D z, $$
where $\Xi = (\infty,\delta] e^{\ii3\pi/2} \cup \delta e^{\ii[3\pi/2, -3\pi/2]} \cup 
[\delta,\infty) e^{-\ii3\pi/2}$ for some sufficiently small $\delta>0$. 
Hence for each $y>0$ the $\calB\big (L_p(\R^{n-1};E)\big)$-valued function $\theta \mapsto L_{\theta} e^{- y L_{\theta}^{1/2m}}$ 
is smooth and all of its derivatives are bounded. So we can apply the inverse Fourier 
transform with respect to $t$ and obtain that
$$ L_{\theta} e^{- y L_{\theta}^{1/2m}} =   \calF_t L e^{- y L^{1/2m}} \calF_t^{-1}$$
on Schwartz functions, where  $L:= 1+ \partial_t + (-\Delta_{n-1})^m$ and  $\calE := e^{-\cdot L^{1/2m}}.$
Here, for $X\in \{E,F\}$  we consider  $L$ as an operator on 
$L_{p,\mu}\big( \R_+; L_p(\R^{n-1};X)\big)$ with the domain 
$$D(L) = D(\partial_t)+D((-\Delta_{n-1})^m)= {}_0W_{p,\mu}^1\big(\R_+; L_p(\R^{n-1};X)\big) 
            \cap L_{p,\mu}\big(\R_+; W_p^{2m}(\R^{n-1};X)\big).$$
In Lemma~3.1 of \cite{MS11} we have established that $L$ is  invertible and sectorial with 
angle not larger than $\pi/2$. Similarly one can treat fractional powers and derive
$L_\theta^{-m_j/2m}= \calF_tL^{-m_j/2m} \calF_t^{-1}$. We arrive at
$$u = \calL_u g := \big (\calF_{t}^{-1} \calT(\theta)\calF_t\big) \, L \,\calE  \, \big ( \calF_t^{-1}\calF_{x'}^{-1} M_u^0(b,\zeta,\vartheta)\calF_{x'} \calF_t \big) \,\big( L^{-m_j/2m} g_j\big)_{j=0,...,m},$$
Analogous arguments show that the second component $\rho$ can be represented by 
$$\rho = \calL_\rho g :=   S^{-1}  L^{m_0/2m}\, \big ( \calF_t^{-1} \calF_{x'}^{-1} M_{\rho}^0(b,\zeta,\vartheta)\calF_{x'} \calF_t \big ) \big( L^{-m_j/2m} g_j\big)_{j=0,...,m},$$
with the operator
\begin{align*}
S &:= 1+\partial_t + (-\Delta_{n-1})^{l/2}\qquad \text{ in the Cases 1 and 2},\\ 
S &:= 1+\partial_t + \sum_{j\in \calJ} (-\Delta_{n-1})^{k_j/2} L^{(m_0-m_j)/2m}\qquad \text{ in  Case 3}.
\end{align*}
Using the properties of $L$ proved in Lemma 3.1 of \cite{MS11}, it can be shown 
as in Section~4.2 of \cite{DPZ08} that $S$ is an isomorphism  between ${}_0\E_{\rho,\mu}$ and ${}_0\F_{0,\mu}$. 
Because $\calD$ is a dense subset of ${}_0\F_\mu$, it now remains to prove the estimate 
\beq\label{lg9}
|\calL_u g|_{\E_{u,\mu}} + |\calL_\rho g|_{\E_{\rho,\mu}} \lesssim |g|_{{}_0\F_{\mu}}, 
 \qquad g\in \calD.
\eeq 
If \eqref{lg9} has been verified then the solution operator $\calL:= (\calL_u,\calL_\rho)$ 
extends continuously to an operator from ${}_0\F_{\mu}$ to ${}_0\E_{u,\mu}\times{}_0\E_{\rho,\mu}$,
and this extension yields the solution of (\ref{dynmaxregHS}).
\textbf{(IV)}  Lemma~3.1 of \cite{MS11} says that for $s\in (0,1]$ we have
$$
D_{L}(s,p) = {}_0W_{p,\mu}^s\big(\R_+; L_p(\R^{n-1};X)\big) \cap L_{p,\mu}\big(\R_+; W_p^{2ms}(\R^{n-1};X)\big).
$$ 
Therefore, for $j=1,...,m$ the operator $L^{-m_j/2m}$ maps the space 
${}_0\F_{j,\mu}=D_{L}(\kappa_j,p)$  continuously into 
\begin{align*}
{}_0\Y_E := D_L(1-1/2mp,p)={}_0W_{p,\mu}^{1-1/2mp}\big(\R_+; L_p(\R^{n-1};E)\big) 
      \cap L_{p,\mu}\big(\R_+; W_p^{2m-1/p}(\R^{n-1};E)\big).
\end{align*}
The same arguments yield that $L^{-m_0/2m}$ maps ${}_0\F_{0,\mu}$ continuously into ${}_0\Y_F$, 
which is defined as ${}_0\Y_E$ with $E$ replaced by $F$. We next prove that the operator
$$\calM^0 := \calF_t^{-1}\calF_{x'}^{-1} M^0(b,\zeta,\vartheta)\calF_{x'} \calF_t$$
on $\calD$ with the symbol $M^0: D_b\times D_\zeta\times \Sigma \ra 
\calB(F\times E^m, E^{2m}\times F)$ given by 
$$ M^0(b,\zeta,\vartheta):=  \big ( M_u^0(b,\zeta,\vartheta), M_\rho^0(b,\zeta,\vartheta) \big),$$
extends continuously to an element of $\calB \big ( {}_0\Y_F \times {}_0\Y_E^m, 
{}_0\Y_E^{2m}\times {}_0\Y_F\big ).$ To this end, we consider the approximating operators
$$\calM^{0,\eps} := \calF_t^{-1}\calF_{x'}^{-1}  M^{0}(b,\zeta,\vartheta)  (1+ \vartheta)^{-\eps}\calF_{x'} \calF_t, \qquad \eps \in (0,1).$$ 
Observe that $\calM^{0,\eps}(1+L^{1/2m})^\eps=\calM^0$ on $\calD$. 
Cauchy's formula yields the representation
$$\calM^{0,\eps} = -\frac{1}{4\pi^2}\int_{\Xi_{\vartheta}} \int_{\Xi_{b}} \calF_t^{-1}\calF_{x'}^{-1}  M^{0}(\wt{b},\zeta,\wt{\vartheta})  (1+ \wt{\vartheta})^{-\eps} (\wt{b}-b)^{-1} (\wt{\vartheta} -\vartheta)^{-1}\calF_{x'} \calF_t \D \wt{b} \D \wt{\vartheta},$$
with $\Xi_{\vartheta} = (-\infty,0] e^{-\ii\phi_*}\cup [0,\infty) e^{\ii\phi_*}$ for some
 $\phi_*\in (\pi/4m,\phi)$, and where  $\Xi_b$ is a closed curve in $D_b$ surrounding 
 $(\ol{B}_{1/2}(1/2))^{1/2m}$. Since $\zeta = \frac{\xi'}{|\xi'|}$ is independent of 
$\theta$, we may rewrite the above equality as
$$\calM^{0,\eps} = -\frac{1}{4\pi^2}\int_{\Xi_{\vartheta}} \int_{\Xi_{b}} \calF_{x'}^{-1} M^{0}(\wt{b},\zeta,\wt{\vartheta})\calF_{x'} (1+ \wt{\vartheta})^{-\eps} (\wt{b}-B)^{-1} (\wt{\vartheta} -L^{1/2m})^{-1}  \D \wt{b} \D \wt{\vartheta},$$
where $B:= (-\Delta_{n-1})^{1/2}L^{-1/2m}$  corresponds to the symbol 
$b =\frac{|\xi'|}{\vartheta}$. The realization of $B$ on $L_{p,\mu}\big(\R_+;L_p(\R^{n-1};E)\big)$ 
is a bounded operator, and its spectrum is contained in the set $(\ol{B}_{1/2}(1/2))^{1/2m}$. 
This can be seen using the joint functional calculus for $\partial_t$ and $(-\Delta_{n-1})^m$ 
on $L_{p,\mu}\big(\R_+;L_p(\R^{n-1};E)\big)$, see Theorem~4.5 of \cite{KW01}.

Due to the $\calR$-boundedness of the sets (\ref{lg7}) and (\ref{kv38}),
the operator-valued Fourier-multiplier theorem  in $\R^{n-1}$
(Theorem~3.25 of \cite{DHP03}, see also Theorem~4.13 of \cite{KW04})
and real interpolation imply that that the operators 
$$M^1(\wt{b}, \wt{\vartheta}) := \calF_{x'}^{-1} M^{0}(\wt{b},\cdot ,\wt{\vartheta})\calF_{x'}, 
     \qquad  \qquad \wt{b} \in D_b,\qquad \wt{\vartheta}\in \Sigma,$$
extend continuously to elements of $\calB\big (W_p^s(\R^{n-1}; F\times E^{m}), 
W_p^s(\R^{n-1}; E^{2m}\times F)\big)$, $s\geq 0$, with uniformly bounded operators norms.  
Since $M^0$ is holomorphic, also $M^1$ is holomorphic in its arguments. By canonical pointwise 
extension we thus obtain that
 $$M^1:D_b\times \Sigma \ra \calB \big ( {}_0\Y_F \times {}_0\Y_E^m, {}_0\Y_E^{2m}
            \times {}_0\Y_F\big )$$ 
is bounded and holomorphic. Using $L$ as an isomorphism 
$D(L) \ra L_{p,\mu}\big (\R_+; L_p(\R^{n-1};E)\big)$ that commutes with $B$, we see that the
 spectrum of the realization of $B$ on $D(L)$ is also contained in  
$(\ol{B}_{1/2}(1/2))^{1/2m}$. By interpolation, the same holds on  ${}_0\Y_F \times {}_0\Y_E^m$.
Hence, we may rewrite $\calM^{0,\eps}$ as
$$\calM^{0,\eps} = -\frac{1}{4\pi^2} \int_{\Xi_{\vartheta}} \int_{\Xi_{b}} 
     M^1(\wt{b}, \wt{\vartheta}) (1+ \wt{\vartheta})^{-\eps} (\wt{b}-B)^{-1} 
        (\wt{\vartheta} -L^{1/2m})^{-1}  \D \wt{b} \D \wt{\vartheta},$$
where the curve integrals are now defined in $\calB \big ( {}_0\Y_F \times {}_0\Y_E^m, 
{}_0\Y_E^{2m} \times {}_0\Y_F\big ).$ 
We thus obtain
$$\calM^{0,\eps} = \frac{1}{2\pi\ii}\int_{\Xi_{\vartheta}} M^2(\wt{\vartheta}) 
 (1+ \wt{\vartheta})^{-\eps} (\wt{\vartheta} -L^{1/2m})^{-1}   \D \wt{\vartheta}$$ 
for a bounded  holomorphic map 
$$M^2: \Sigma \ra \calB \big ( {}_0\Y_F \times {}_0\Y_E^m, {}_0\Y_E^{2m}\times {}_0\Y_F\big ).$$
Since the realization of $L^{1/2m}$ on $L_{p,\mu}\big (\R_+; L_p(\R^{n-1};E)\big )$ is sectorial with angle not larger than $\pi/4m$, it follows from Corollary~1 of \cite{CP03} that $L^{1/2m}$ admits a  bounded operator-valued $\calH^\infty$-calculus with $\calH^\infty$-angle not larger than $\pi/4m$ on the real interpolation spaces ${}_0\Y_E^m$ and ${}_0\Y_F$, respectively. From this fact and the boundedness of $M^2$ on $\Sigma$ we infer
\beq\label{kv39}
|\calM^{0,\eps}|_{\calB ( {}_0\Y_F \times {}_0\Y_E^m, {}_0\Y_E^{2m}\times {}_0\Y_F)} \lesssim \sup_{\wt{\vartheta}\in \Sigma} |M^2(\wt{\vartheta}) (1+ \wt{\vartheta})^{-\eps}|_{\calB  ( {}_0\Y_F \times {}_0\Y_E^m, {}_0\Y_E^{2m}\times {}_0\Y_F )} \leq C,  
\eeq 
where $C$ does not depend on $\eps\in (0,1)$. Due to Proposition~2.2 of \cite{DHP03}, for $h\in D(L^2)$ the map $\eps \mapsto (1+L^{1/2m})^\eps h$ is continuous with values in $D_{L}(1-1/2mp,p)$. Together with (\ref{kv39}), this fact yields 
\begin{align*}
|\calM^0 h|_{{}_0\Y_E^{2m}\times {}_0\Y_F}&\,  \lesssim \limsup_{\eps\ra 0} |\calM^{0,\eps}|_{\calB ( {}_0\Y_F \times {}_0\Y_E^m, {}_0\Y_E^{2m}\times {}_0\Y_F)} \,|(1+L^{1/2m})^\eps h|_{{}_0\Y_F \times {}_0\Y_E^m}  \lesssim |h|_{{}_0\Y_F \times {}_0\Y_E^m}.
\end{align*}
Since $D(L^2)$ is dense in $D_{L}(1-1/2mp,p)$, we obtain that $\calM^0$ extends to an element of $\calB \big ( {}_0\Y_F \times {}_0\Y_E^m, {}_0\Y_E^{2m}\times {}_0\Y_F\big )$, as asserted.

\textbf{(V)} Now we can show the required estimate for $\calL_u$, i.e.,
\beq\label{kv7}
|\calL_u g|_{\E_{u,\mu}}\lesssim |g|_{{}_0\F_\mu}, \qquad g\in \calD.
\eeq
The extension operator $\calE = e^{-\cdot L^{1/2m}}$ maps continuously 
$$D_L(1-1/2mp,p) = D_{L^{1/2m}}(2m-1/p,p) \ra L_p\big(\R_+; D(L)\big),$$ 
and $L$ maps the space $L_p\big(\R_+; D(L)\big)$ continuously into 
$$L_p\big (\R_+; L_{p,\mu}(\R_+; L_p(\R^{n-1};E))\big ) = L_{p,\mu}\big(\R_+; L_p(\R_+^n;E)\big).$$ 
Of course, here $E$ may be replaced by $F$. Thus $L\,\calE$ maps continuously
$${}_0\Y_E^{2m}\times {}_0\Y_F \ra L_{p,\mu}\big(\R_+; L_p(\R_+^n;E^{2m}\times F)\big).$$
Theorem~\ref{sec:mumult} and the $\calR$-boundedness of (\ref{lg6}) imply that
 $\calF_{t}^{-1} \calT(\cdot) \calF_t$ extends to a continuous operator
$$L_{p,\mu}\big (\R_+; L_p(\R_+^n;E^{2m})\big) \ra L_{p,\mu}\big(\R_+; W_p^{2m}(\R_+^n;E)\big).$$ 
Alternatively, this fact follows from Proposition~\ref{sec:mumultprop} since
 one can show that the operator family 
$$\big \{\theta^2 \frac{\partial}{\partial\theta^2}  D^\alpha \calT(\theta)\;:\; \theta\in \R,\;\; |\alpha|\leq 2m\big \}$$ 
is bounded in $\calB\big( L_p(\R_+^n;E^{2m}), L_p(\R_+^n;E)\big)$ arguing 
as in the proof of Lemma~4.4 of \cite{DHP07}.
The equation for $u$ shows that its $\E_{u,\mu}$-norm  can be controlled by its 
$L_{p,\mu}\big(\R_+; W_p^{2m}(\R_+^n;E)\big)$-norm. So we have established (\ref{kv7}). We finally 
consider the required estimate for $\calL_\rho$. As above we obtain that $L^{m_0/2m}$ maps 
continuously 
$${}_0 \Y_F = D_{L}(1-1/p,p) \ra D_{L}(\kappa_0,p) = {}_0\F_{0,\mu}.$$  
Since $S^{-1}$ is an isomorphism from ${}_0\F_{0,\mu}$ to ${}_0 \E_{\rho,\mu}$, this gives the estimate
for $\calL_\rho$. \eprf

The analogous half-space result for (\ref{statmrI}) reads as follows.

\begin{lem}\label{sec:HSP2}\textsl{Let $E$ be a Banach space of class $\HT$, $p\in (1,\infty)$, $\mu\in (1/p,1]$, and assume that \emph{(E)} and \emph{(LS)} are valid. Then for $(f,\wt{g},u_0) \in \calD_{\emph{\text{stat}}}$ there is a unique solution $u\in \E_{u,\mu}$ of 
\begin{alignat}{3}
u +  \partial_t u + \calA(D) u & =    f(t,x), & \qquad &  x\in \R_+^n, &\qquad &  t>0,  \nonumber\\
\calB_j(D) u & =  g_j(t,x), && x\in \R^{n-1}, &&  t>0,\qquad j=1,...,m, \label{eq:HSP2}\\
u(0,x) & =  u_0(x), \quad  && x\in  \R_+^n. && \tag*{\BlackBox}
\end{alignat}
}
\end{lem}

We refrain from giving a detailed proof of this result, which is similar to the one 
of Lemma~\ref{lg10}  and also less sophisticated. (See Section 2 of \cite{Mey10} for the details.)
Again we may restrict to the case $f=0$, 
$\wt{g}\in {}_0\wt{\F}_\mu$ and $u_0=0$. Applying the partial Fourier transforms with respect to $t$ and $x'$ to 
(\ref{eq:HSP2}) we arrive at an ordinary initial value problem, whose solution operator 
is for regular  data $(g_1,..., g_m)$ given by 
$$\wt{\calL} = \wt{\calT} \big( L^{1-m_j/2m} \calE g_j\big)_{j=1,...,m},$$ 
due to Lemma~4.3 of \cite{DHP07}. Here $\wt{\calT}$ has the same properties as $\calT$ and 
$L$, $\calE$ are given as in the proof of Lemma~\ref{lg10}. The arguments given in the 
Steps~IV and V of the proof above yield that $\wt{\calL}\in \calB( {}_0\wt{\F}_{\mu}, 
{}_0\E_{u,\mu})$, which implies the solvability of (\ref{eq:HSP2}) as asserted.

\section{The general problem on a domain}
Theorems \ref{statmrthm} and \ref{dynmrthm} are now a consequence of the above results for the 
model problems and a perturbation and localization procedure, analogous to the one in e.g.\
Section~4.5 of \cite{DPZ08}. We only sketch the proof below since the full procedure
is rather lengthy and tedious. The arguments are worked out in great detail 
in Sections 2.3, 2.4, 3.2.2 and 3.3 of \cite{Mey10}.
Moreover, we concentrate on (\ref{dynmaxregI}) since the proof for (\ref{statmrI}) is
similar and a bit simpler.

\textbf{Proof of Theorem \ref{dynmrthm}.} 
\textbf{(I)} Let us first consider the necessary conditions on the data. The considerations in 
Section~2 and the assumptions (SD), (SB) and (SC) yield that $\calA\in \calB\big( \E_{u,\mu}, 
\E_{0,\mu}\big)$ and $\calB_j \in \calB\big( \E_{u,\mu}, \F_{j,\mu}\big)$, $\calC_j \in 
\calB\big( \E_{\rho,\mu}, \F_{j,\mu}\big)$ for $j=0,...,m$. Moreover, we have 
$$W_{p,\mu}^{\kappa_j}\big (J; L_p(\Gamma;E)\big)\hra BU\!C\big (\ol{J}; L_p(\Gamma;E)\big) 
 \qquad \text{if }\; \kappa_j > 1-\mu+1/p$$ 
 for $j=1,...,m$, due to Proposition~2.10 of \cite{MS11}. Thus in this case the $j$-th boundary 
equation in (\ref{dynmaxregI}) must hold up to $t=0$ by continuity, which explains the 
compatibility conditions in $\calD_{\text{rel}}$ for this case. Similiarly, for $j=0$ the
regularity compatibility at the boundary is needed if $\kappa_0 > 1-\mu+1/p$, i.e., if 
$\partial_t \rho$ has a trace at $t=0$. For the existence of a solution 
$(u,\rho)\in \E_{u,\mu}\times \E_{\rho,\mu}$, it is therefore necessary that the data in 
(\ref{dynmaxregI}) belong to $\calD_{\text{rel}}$.

\textbf{(II)} Let us show that $(f,g,u_0,\rho_0)\in \calD_{\text{rel}}$ is also sufficient for the 
existence of a unique solution $(u,\rho)\in \E_{u,\mu}\times \E_{\rho,\mu}$ of (\ref{dynmaxregI}). 
Uniqueness follows as in Step I of the proof of Lemma \ref{lg10}. For the existence of the solution
 $(u,\rho)$, note that it suffices to consider small $T>0$ by a standard compactness argument.
For simplicity, we assume that $\Omega$ is bounded. The case of unbounded $\Omega$ requires 
minor modifications.

We cover $\ol{\Omega}$ by a finite number of open balls $B_i$ such that 
$B_i \cap \Gamma =  \emptyset$ for $i=1,..., N_F$ and $B_i \cap \Gamma \neq \emptyset$ for 
$i=N_H+1,..., N_F$, where $N_F, N_H\in \N$. We further take a smooth partition of unity $\psi_i$ 
for $\ol{\Omega}$ subordinate to this cover. Let  $(u,\rho)\in \E_{u,\mu}\times \E_{\rho,\mu}$.
Now, $(u,\rho)$ solves (\ref{dynmaxregI}) if and only if $(u_i,\rho_i) = (\psi_i u, \psi_i \rho)$ 
satisfies 
\begin{alignat}{3}
\partial_t u_i + \calA u_i & =    \psi_i f + [\calA, \psi_i] u, & \quad & \text{in }\; \Omega\cap B_i, & \quad &  t\in J,  \nonumber\\
\partial_t \rho_i + \calB_0 u_i + \calC_0\rho_i & =  \psi_i g_0 + [\calB_0, \psi_i] u  + [\calC_0,\psi_i]\rho, &&\text{on }\; \Gamma \cap B_i,  &&  t\in J,\quad   \nonumber \\
\calB_ju_i + \calC_j\rho_i & =  \psi_i g_j + [\calB_0, \psi_i] u  + [\calC_0,\psi_i]\rho, &&\text{on }\; \Gamma \cap B_i,  &&  t\in J,\quad j=1,...,m, \label{loc1}\\
u_i|_{t=0} & =  \psi_i u_0, &&  \text{in }\; \Omega\cap B_i, &&  \nonumber \\
\rho_i|_{t=0} & =  \psi_i\rho_0, &&  \text{on }\; \Gamma \cap B_i, &&  \nonumber
\end{alignat}
for all $i=1,..., N_H$, where  $[\calA, \psi_i] u=\calA (\psi_i u) - \psi_i \calA u$.
For $i=1,..., N_F$ no boundary conditions are involved in (\ref{loc1}). We extend the coefficients 
of $\calA$ outside $B_i$ to $\R^n$ such that (SD) is still valid, and denote the operator with 
extended coefficients by $\calA^i$. Then $u_i$ solves (\ref{loc1}) for $i=1,..., N_F$ if and only 
if it solves
\begin{alignat}{3}
\partial_tu_i + \calA^i u_i & =  \psi_i f + [\calA, \psi_i] u, & \qquad & \text{in }\; \R^n, \qquad t\in J, \nonumber\\
u_i|_{t=0} & =  \psi_i u_0, && \text{in }\; \R^n.   \label{loc2}
\end{alignat}
Due to the continuity of the top order coefficients of $\calA$, the top order
part of the operator $\calA^i$ is a small perturbation of a homogeneous constant coefficient operator 
satisfying (E) if the extension of the coefficients 
is appropriate, provided $T$ and the radius of $B_i$ are sufficiently small.  
Poincar\'e's inequality in the $L_{p,\mu}$-spaces (Lemma~2.12 in \cite{MS11}) allows to estimate 
lower order terms
with constants decreasing to 0 as $T\to 0$, see Lemma~1.3.13 of \cite{Mey10}.
Using Lemma~\ref{sec:HSP2easy}, we can now solve (\ref{loc2}) by a straightforward fixed point 
argument. We thus obtain a continuous solution operator 
$\calL_F^i:\E_{0,\mu}(J\times \R^n)\times X_{u,\mu}(\R^n)\ra \E_{u,\mu}(J\times \R^n)$ 
for (\ref{loc2}). It follows that
$$u_i = \calL_F^i\big ( \psi_i f + [\calA, \psi_i] u, \psi_i u_0\big ), \qquad i=1,..., N_F.$$
Observe that the commutator terms are of lower order.
For $i=N_F+1,..., N_H$ the boundary conditions in (\ref{loc1}) are present. We choose the 
$B_i$ so small that we have a chart $\vphi_i$ for $\Gamma$ with domain $B_i$ associated to some 
$x_i \in \Gamma$. Denoting by $\Phi_i$ the corresponding push-forward operator, i.e., 
$\Phi_i v = v\circ  \vphi_i^{-1}$, we obtain that $(u_i,\rho_i)$ solves  (\ref{loc1}) 
if and only if $(v_i, \sigma_i ) = (\Phi_i u_i, \Phi_i \rho_i)$ solves
\begin{alignat}{3}
\partial_t v_i + \big ( \Phi_i \calA \Phi_i^{-1}\big)   v_i & =    \Phi_i\big( \psi_i f + [\calA, \psi_i] u \big), & \quad & \text{in }\; \R_+^n\cap \vphi_i(B_i), & \quad &  \nonumber\\
\partial_t \sigma_i + \big ( \Phi_i \calB_0 \Phi_i^{-1}\big)  v_i + \calC_0^{\g_i}  \sigma_i & =   \Phi_i\big( \psi_i g_0 + [\calB_0, \psi_i] u  + [\calC_0,\psi_i]\rho\big), &&\text{on }\; \R^{n-1}\cap \vphi_i(B_i),  &&     \nonumber \\
\big ( \Phi_i \calB_j \Phi_i^{-1}\big)  v_i+ \calC_j^{\g_i} \sigma_i & =  \Phi_i\big(\psi_i g_j + [\calB_j, \psi_i] u  + [\calC_j,\psi_i]\rho\big), &&\text{on }\; \R^{n-1}\cap \vphi_i(B_i),  &&  \nonumber\\
v_i|_{t=0} & =  \Phi_i \psi_i u_0, &&  \text{in }\; \R_+^n\cap \vphi_i(B_i), &&  \nonumber \\
\sigma_i|_{t=0} & =  \Phi_i \psi_i \rho_0, &&  \text{on }\; \R^{n-1}\cap \vphi_i(B_i), &&  \nonumber
\end{alignat}
for $t\in J$ and $j=1,...,m$.  Recall that $\calC_j^{\g_i}$ denotes the local representation of 
$\calC_j$ with respect to the coordinates $\g_i$ corresponding to $\vphi_i$. According to 
Theorem~10.3 of \cite{Wlo87}, at $t\in \ol{J}$ and $x_i$ the principal parts of the operators 
$\Phi_i \calA \Phi_i^{-1}$ and $\Phi_i \calB_j \Phi_i^{-1}$ are given by
$$\calA_\sharp\big(t, x_i, \calO_{\nu(x_i)}^T D\big),\qquad 
  \calB_{j\sharp}\big(t, x_i, \calO_{\nu(x_i)}^TD\big),$$ 
respectively. Extending now the coefficients of the transformed operators 
$\Phi_i \calA \Phi_i^{-1}$, $\Phi_i \calB_j \Phi_i^{-1}$ and  $\calC_j^{\g_i}$ such that 
(SD), (SB) and (SC) remain valid, we obtain that $(\Phi_i u_i, \Phi_i \rho_i)$ solves a half-space 
problem with operators that are either of lower order or small perturbations of 
constant coefficient operators satisfying the conditions of Lemma~\ref{lg10}. As for the 
full-space case, if $T$ and $B_i$ are sufficiently small, then a continuous solution operator 
$\calL_H^i$ exists for this half-space problem, which maps the relevant data space 
continuously into  $\E_{\rho,\mu}(J\times \R^n_+)\times  \E_{\rho,\mu}(J\times \R^n_+)$. 
For $i=N_F+1, ..., N_H$ we thus  obtain
$$(u_i, \rho_i) = \Phi_i^{-1}R_i \calL_H^i\big (\Phi_i( \psi_i f + [\calA, \psi_i] u), 
\Phi_i( \psi_i g + [\calB, \psi_i] u  + [\calC,\psi_i]\rho), \Phi_i \psi_i u_0, 
\Phi_i \psi_i \rho_0 \big).$$
Here $R_i$ is the restriction to $\R^n_+\cap\vphi_i(B_i)$ and 
we have set $\calB= (\calB_0,..., \calB_m)$ and $\calC= (\calC_0,..., \calC_m)$ 
for simplicity. Again, the commutator terms are of lower order.

\textbf{(III)} We next choose smooth functions $\phi_i$, $i=1,..., N_H$, satisfying 
$\phi_i \equiv 1$ on $\supp\,\psi_i$ and $\supp\, \phi_i \subset B_i$. The above considerations 
show that if $(u,\rho)$ solves  (\ref{dynmaxregI}) then it is a fixed point of the map $\calG_{f,g,u_0,\rho_0}(u,\rho) := \sum_i \phi_i (u_i,\rho_i)$ on the complete metric space
$$Z_{u_0,\rho_0} := \big\{ (u,\rho)\in \E_{u,\mu}\times \E_{\rho,\mu}\;:\; u|_{t=0}
 = u_0, \;\;\rho|_{t=0}= \rho_0\big\}.$$
We remark that $Z_{u_0,\rho_0}$ is nonempty by Lemma~4.4 of \cite{MS11} and Lemma 3.2.2 
of \cite{Mey10}. Since the operators in the arguments of $\calL_F^i$ and $\calL_H^i$ are of lower 
order, one can show that for all data $(f^*,g^*,u_0^*,\rho_0^*)\in \calD_{\text{rel}}$ 
the map $\calG_{f^*,g^*,u_0^*,\rho_0^*}$ has indeed a unique fixed point on 
$Z_{u_0^*,\rho_0^*}$, making $T$ and $B_i$ once more smaller if necessary. Another fixed point 
argument yields for given data $(f,g,u_0,\rho_0)\in \calD_{\text{rel}}$ the appropriate auxiliary 
data $(f^*,g^*,u_0^*,\rho_0^*)\in \calD_{\text{rel}}$ such that the fixed point of 
$\calG_{f^*,g^*,u_0^*,\rho_0^*}$ is the solution of (\ref{dynmaxregI}).

\textbf{(IV)} To finish the proof, note that the continuity of the resulting solution operator 
$\calL_{\text{rel}}$ for (\ref{dynmaxregI}) is a consequence of the open mapping theorem. 
Moreover, the norm of $\calL_{\text{rel}}$ restricted to $\calD_{\text{rel}}^0$ 
is uniform in $T$ due to an extension argument. It uses the extension operator from 
Lemma~2.5 of \cite{MS11} for the ${}_0W_{p,\mu}^s$-spaces over $J$ to the half-line, 
whose norm is independent of the length of $J$. \hfill \BlackBox

\end{document}